\documentclass[a4paper,11pt]{amsart}
\usepackage[francais]{babel}
\usepackage{epsfig,xypic}
\usepackage{amssymb,bbm}
\theoremstyle{plain}
\author{Ivan Marin}
\address{69 rue S\'ebastien Gryphe\\ F-69007 Lyon}
\email{imarin@maths.univ-evry.fr}
\urladdr{http://www.maths.univ-evry.fr/pages\_perso/marin/}
\title[Caract\`eres de rigidit\'e de GT]{Caract\`eres de rigidit\'e du groupe de Grothendieck-Teichm\"uller}

\newtheorem*{THEOA}{Th\'eor\`eme A}
\newtheorem*{THEOB}{Th\'eor\`eme B}
\newtheorem*{THEOC}{Th\'eor\`eme C}
\newtheorem*{THEOD}{Th\'eor\`eme D}
\newtheorem{lemme}{Lemme}
\newtheorem{defi}{D\'efinition}
\newtheorem{prop}{Proposition}
\newtheorem*{ccor}{Corollaire}

\newtheorem{theo}{ThŽorme}

\def\stvrule{\vrule width-0.4pt}         
\def\sthrule{\hrule \hrule height-0.4pt} 

\def\ststrut{\vrule height1.6ex 
                     depth0.6ex 
                       width0ex 
\relax}
\def\scarre#1{%
    \vcenter{\hbox{}\hrule
             \hbox{\vrule\makebox[2.3ex]{\ststrut$\scriptstyle#1$}\vrule}\sthrule}%
             \stvrule}
\def\sgenruban#1{\vcenter{\halign{&$\scarre{##}$\cr#1}}\egroup}

\def\smallruban{%
  \bgroup
  \let\ =\omit
  \let\\=\cr
  \offinterlineskip
  \sgenruban}

\def\Stvrule{\vrule width-0.4pt}         
\def\Sthrule{\hrule \hrule height-1.6pt} 

\def\Ststrut{\vrule height0.9ex 
                     depth0.3ex 
                       width0ex 
\relax}
\def\Scarre#1{%
    \vcenter{\hbox{}\hrule
             \hbox{\vrule\makebox[1.0ex]{\Ststrut$\scriptstyle#1$}\vrule}\Sthrule}%
             \Stvrule}

\def\Sgenruban#1{\vcenter{\halign{&$\Scarre{##}$\cr#1}}\egroup}

\def\Smallruban{%
  \bgroup
  \let\ =\omit
  \let\\=\cr
  \offinterlineskip
  \Sgenruban}

\def\rtvrule{\vrule width-0.8pt}         
\def\rthrule{\hrule \hrule height-0.8pt} 

\def\rtstrut{\vrule height3.2ex 
                     depth1.2ex 
                       width0ex 
\relax}
\def\rcarre#1{%
    \vcenter{\hbox{}\hrule
             \hbox{\vrule\makebox[4.6ex]{\rtstrut$\scriptstyle#1$}\vrule}\rthrule}%
             \rtvrule}
\def\rgenruban#1{\vcenter{\halign{&$\rcarre{##}$\cr#1}}\egroup}

\def\ruban{%
  \bgroup
  \let\ =\omit
  \let\\=\cr
  \offinterlineskip
  \rgenruban}

\newcommand{\la}{\lambda}

\def\Hom{\mathrm{Hom}}

\def\Ker{\mathrm{Ker}}

\def\hh{\hslash}

\def\SN{\ensuremath{\mathfrak{S}_n}}
\def\C{\ensuremath{\mathbbm{C}}}
\def\Q{\mathbbm{Q}}
\def\Z{\mathbbm{Z}}
\def\R{\mathbbm{R}}
\def\N{\mathbbm{N}}

\def\ot{\otimes}
\def\k{\mathbbm{k}}
\def\Ass{\mathbbm{Ass}}
\def\AN{\ensuremath{\mathfrak{B}_n}}

\def\eps{\epsilon}

\def\GM{\mathbbm{G}_m}

\def\into{\hookrightarrow}
\def\ii{\mathrm{i}}
\def\kt{\k^{\times}}
\def\bar{\overline}
\def\Id{\mathrm{Id}}
\def\cc{\mathfrak{c}}
\def\Gal{\mathrm{Gal}}
\def\lmm{ \left( \begin{array}{cc}}
\def\rmm{\end{array} \right)}
\def\diag{\mathrm{diag}}
\begin{document}

\maketitle

\noindent {\bf R\'esum\'e.} Soit $\k$ un corps (topologique) de caract\'eristique
nulle. \`A l'aide d'un associateur de Drinfeld $\Phi$
on peut associer \`a toute re\-pr\'e\-sen\-ta\-tion $\rho$ d'une certaine
$\k$-alg\`ebre de Hopf $\mathfrak{B}_n(\k)$ une re\-pr\'e\-sen\-ta\-tion
$\widehat{\Phi}(\rho)$ du groupe de tresses sur le corps $\k((h))$
des s\'eries de Laurent. Nous \'etudions la d\'e\-pen\-dance en
$\Phi$ de $\widehat{\Phi}(\rho)$ pour certaines re\-pr\'e\-sen\-ta\-tions,
dites GT-rigides, et en d\'eduisons des re\-pr\'e\-sen\-ta\-tions projectives
(continues) du groupe de Grothendieck-Teich\-m\"uller $GT_1(\k)$,
donc pour $\k = \Q_l$ du groupe de Galois absolu de $\Q(\mu_{l^{\infty}})$. La
plupart du temps, ces re\-pr\'e\-sen\-ta\-tions projectives se d\'ecomposent
en caract\`eres lin\'eaires, que nous d\'eterminons pour les
re\-pr\'e\-sen\-ta\-tions de l'alg\`ebre d'Iwahori-Hecke de type A. Dans ce
cas, nous calculons de plus $\widehat{\Phi}(\rho)$ pour $\Phi$
un associateur pair, et en d\'eduisons un mod\`ele
matriciel unitaire des re\-pr\'e\-sen\-ta\-tions de l'alg\`ebre d'Iwahori-Hecke.
Pour l'action de $GT_1(\k)$, les re\-pr\'e\-sen\-ta\-tions de cette alg\`ebre
qui correspondent aux diagrammes de Young \og en \'equerre \fg\ jouissent
de propri\'et\'es remarquables.

\bigskip

\noindent {\bf Abstract.} Let $\k$ be a (topological) field of characteristic 0. Using a Drinfeld associator, a representation $\widehat{\Phi}(\rho)$
of the braid group over the field $\k((h))$ of Laurent series can be
associated to any representation of a certain Hopf algebra $\mathfrak{B}_n(\k)$.
We investigate the dependance in $\Phi$ of $\widehat{\Phi}(\rho)$ for
a certain class of representations --- so-called GT-rigid representations ---
and deduce from it (continuous) projective representations of the
Grothendieck-Teichm\"uller group $GT_1(\k)$, hence for $\k = \Q_l$
representations of the absolute Galois group of $\Q(\mu_{l^{\infty}})$.
In most situations, these projective representations can be decomposed
into linear characters, which we do for the representations of the
Iwahori-Hecke algebra of type A. In this case, we moreover express
$\widehat{\Phi}(\rho)$ when $\Phi$ is even, and get unitary
matrix models for the representations of the Iwahori-Hecke algebra.
With respect to the action of $GT_1(\k)$, the representations
of this algebra corresponding to hook diagrams have noticeable
properties.

\bigskip


\section{Introduction}

V.G. Drinfeld a introduit dans \cite{DRIN}, pour tout corps $\k$ de
caract\'eristique 0 et tout $\la \in \k$ une famille $\Ass_{\la}(\k)$
de s\'eries formelles $\Phi$ en deux variables non commutatives qui
permettent de construire des morphismes $\widetilde{\Phi}$
d'un compl\'et\'e $B_n(\k)$ du groupe de tresses $B_n$ vers le groupe
des inversibles de la compl\'etion
d'une certaine $\k$-alg\`ebre de Hopf gradu\'ee $\mathfrak{B}_n(\k)$.
Ces s\'eries sont appel\'ees des associateurs (de Drinfeld). On en
d\'eduit dans \cite{ASSOC} des foncteurs $\widehat{\Phi}$ de la cat\'egorie
des repr\'esentations (ici toujours suppos\'ees de dimension finie) de
$\mathfrak{B}_n(\k)$ vers la cat\'egorie des repr\'esentations
de $B_n$ sur le corps des s\'eries de Laurent $K = \k((h))$.
L'int\'er\^et de ces foncteurs est, d'une part, que les repr\'esentations
de $\mathfrak{B}_n(\k)$ sont plus ais\'ees \`a \'etudier et \`a
construire que celles de $B_n$, et d'autre part que ces foncteurs
jouissent de bonnes propri\'et\'es en termes de th\'eorie des
repr\'esentations. Notamment, l'irr\'eductibilit\'e et
l'absolue irr\'eductibilit\'e des repr\'esentations est pr\'eserv\'ee. De
plus, certains associateurs (les associateurs \og r\'eels \fg, c'est-\`a-dire
si $\k
\subset \R$) permettent de construire des repr\'esentations unitaires
de $B_n$ --- sans pour autant donner de mod\`ele matriciel, l'explicitation
de tels associateurs \'etant en g\'en\'eral d\'elicate.

La premi\`ere motivation de ce travail est d'expliciter ce qui
r\'esulte d'un changement d'associateur. Si $\rho$ est une des
repr\'esentations irr\'eductibles usuelles de $\mathfrak{B}_n(\k)$,
les repr\'esentations $\widehat{\Phi}_1(\rho)$ et $\widehat{\Phi}_2(\rho)$
de $B_n$ pour $\Phi_1,\Phi_2 \in \Ass_{\la}(\k)$ et $\la \neq 0$
sont conjugu\'ees par des endomorphismes diagonalisables. Nous
nous int\'eressons ici \`a ces endomorphismes, notamment dans le
cas des repr\'esentations de l'alg\`ebre d'Iwahori-Hecke de type A, quotient
classique de l'alg\`ebre de groupe de $B_n$.

La deuxi\`eme motivation concerne le groupe de Grothendieck-Teichm\"uller
introduit par Drinfeld dans \cite{DRIN}. Ce groupe se d\'ecline
en trois versions, le groupe profini $\widehat{GT}$, sa
composante pro-$l$ not\'ee $GT^{(l)}$, et la version $\k$-pro-unipotente
$GT(\k)$. Ils se d\'ecomposent comme produit semi-direct
d'un tore et d'un sous-groupe $\widehat{GT}_1$, $GT^{(l)}_1$ ou
$GT_1(\k)$. L'int\'er\^et arithm\'etique de ces groupes est que le
groupe de Galois absolu de $\Q$ se plonge naturellement
dans $\widehat{GT}$. On a d'autre
part des morphismes canoniques $\widehat{GT} \to GT^{(l)} \into
GT(\Q_l)$, dont on d\'eduit un morphisme du groupe de Galois absolu
de $\Q(\mu_{l^{\infty}})$ vers $GT_1(\Q_l)$.

Ce groupe $\widehat{GT}$ (resp. $GT(\k)$) s'identifie \`a un groupe
d'automorphismes d'une compl\'etion pro-finie $\widehat{B_n}$
de $B_n$ (resp. de la compl\'etion $B_n(\k)$), de m\^eme que $B_n$
s'identifie, par l'action d'Artin, \`a un groupe d'automorphismes
d'un groupe libre $F$. Dans ce dernier cas, on obtient classiquement
des repr\'esentations projectives de $B_n$ \`a partir de certaines
repr\'esentations de $F$, les \og syst\`emes locaux rigides \fg, 
c'est-\`a-dire les repr\'esentations irr\'eductibles de $F$ dont
la classe d'isomorphisme est inchang\'ee apr\`es torsion par
l'action de $B_n$. De fa\c con analogue, nous associons ici des
repr\'esentations projectives $Q_R$ de $GT_1(\k)$ \`a certaines
repr\'esentations $R = \widehat{\Phi}(\rho)$, qui pr\'esentent
de plus la particularit\'e que $Q_R$ se d\'ecompose en caract\`eres
lin\'eaires \`a valeurs dans $K^{\times}$. Parmi ces
repr\'esentations, dites GT-rigides et agr\'egeantes, se trouvent
notamment les repr\'esentations irr\'eductibles de
l'alg\`ebre d'Iwahori-Hecke de type A. Si $\k = \Q_l$, les caract\`eres obtenus
induisent enfin des caract\`eres lin\'eaires du groupe de Galois absolu
de $\Q(\mu_{l^{\infty}})$ \`a valeurs dans $\Q_l((h))^{\times}$.

\medskip

On note $P_n$ le groupe de tresses pures. Les r\'esultats
principaux de l'article sont les suivants.

\begin{THEOA} Soit $\rho$ une repr\'esentation absolument irr\'eductible,
GT-rigide et agr\'egeante de $\mathfrak{B}_n(\k)$, $R = \widehat{\Phi}(\rho)$
pour un certain $\Phi \in \Ass_{\la}(\k)$, $\la \neq 0$. $Q_R$
est une repr\'esentation projective continue de $GT_1(\k)$,
non triviale si et seulement si $R([P_n,P_n]) = \{ 1 \}$.  Si
$\k = \Q_l$, $Q_R$ induit une repr\'esentation continue du groupe
de Galois absolu de $\Q(\mu_{l^{\infty}})$, non triviale ssi
$R([P_n,P_n]) = \{ 1 \}$.
\end{THEOA}

D'apr\`es \cite{ASSOC}, un moyen d'obtenir des repr\'esentations
unitaires de $B_n$ consiste \`a calculer $\widehat{\Phi}(\rho)$ pour $\rho$
une repr\'esentation de $\mathfrak{B}_n(\k)$ satisfaisant certaines
conditions de sym\'etrie et $\Phi \in \Ass_{\la}(\k)$ un
associateur tel que $\k \subset \R$ et $\la \neq 0$. En
particulier, l'unitarit\'e des repr\'esentations de l'alg\`ebre d'Iwahori-Hecke est
expliqu\'ee par cette construction \`a partir de re\-pr\'e\-sen\-ta\-tions
par\-ti\-cu\-li\`eres de $\mathfrak{B}_n(\k)$, les re\-pr\'e\-sen\-ta\-tions
in\-fi\-ni\-t\'e\-si\-ma\-les de l'al\-g\`e\-bre d'Iwahori-Hecke. En revanche, aucun
associateur r\'eel n'\'etant \`a l'heure actuelle explicitement connu,
il est en g\'en\'eral difficile de calculer $\widehat{\Phi}(\rho)$.
N\'eanmoins, on sait qu'existe $\Phi \in \Ass_1(\Q)$ \emph{pair},
c'est-\`a-dire tel que \\
$\Phi(-x,-y) = \Phi(x,y)$.

\begin{THEOB} Si $\rho$ est une repr\'esentation infinit\'esimale
de l'alg\`ebre d'Iwahori-Hecke, et $\Phi \in \Ass_{\la}(\k)$ pour $\la \neq 0$
est pair, un mod\`ele matriciel de $\widehat{\Phi}(\rho)$ est donn\'e
par la proposition \ref{modunitaire}. En particulier, si $\k \subset \R$,
ces repr\'esentations sont unitaires au sens de \cite{ASSOC}.
\end{THEOB}

L'\'etude de $Q_R$ pour $R$ la repr\'esentation de Burau r\'eduite
(avec $n$ variable) fait appara\^\i tre une famille de caract\`eres
$\chi_d$, $d \geq 2$, de $GT_1(\k)$ vers $(\k[[h]])^{\times}$.

\begin{THEOC} Les caract\`eres $\chi_d$, $d \geq 2$, sont alg\'ebriquement
ind\'ependants. Si $R$ est une repr\'esentation irr\'eductible
de l'alg\`ebre d'Iwahori-Hecke de type A, les caract\`eres intervenant dans
la d\'e\-com\-po\-si\-tion de $Q_R$ sont des mon\^omes (explicitement
d\'etermin\'es) en les $\chi_d$.
\end{THEOC}

Une repr\'esentation irr\'eductible $R$ de l'alg\`ebre d'Iwahori-Hecke associ\'ee
\`a un diagramme de Young $\alpha$ \'etant donn\'ee, une situation
remarquable se produit lorsque les caract\`eres qui interviennent
dans la d\'ecomposition de $Q_R$ sont deux \`a deux distincts. On
dit alors que $Q_R$ est sans r\'esonances.

\begin{THEOD} La repr\'esentation $Q_R$ est sans r\'esonances
si et seulement si $\alpha$ est en \'equerre ou correspond \`a la
partition [2,2].
\end{THEOD}

Cet article comporte trois parties, qui exposent successivement
les pr\'e\-li\-mi\-nai\-res n\'e\-ces\-saires, les pro\-pri\'e\-t\'es g\'e\-n\'erales
des actions de $GT_1(\k)$ obtenues, et enfin l'\'etude particuli\`ere
aux repr\'esentations de l'alg\`ebre d'Iwahori-Hecke. Le
th\'eor\`eme A est d\'emontr\'e en sections 3.4 et 3.5 (corollaire
de la proposition \ref{homeo} et proposition \ref{trivia}). Le
th\'eor\`eme B est d\'emontr\'e en section 4.4 (proposition
\ref{modunitaire}), le th\'eor\`eme C en section 4.5 et le th\'eor\`eme
D en section 4.6 (proposition \ref{equerres}). Certains calculs utiles
sont repouss\'es en appendice, ainsi qu'une g\'en\'eralisation
de ces repr\'esentations projectives de $GT_1(\k)$ en des
1-cocycles de $GT(\k)$.

\section{Pr\'eliminaires}

\subsection{Tresses et tresses infinit\'esimales}

On note $B_n$ pour $n \geq 1$ le groupe de tresses \`a $n$ brins,
en convenant $B_1 = \{ e \}$. On note $P_n$ le groupe de tresses pures,
noyau de la surjection canonique $ \pi : B_n \to \SN$, $\SN$ d\'esignant le 
groupe sym\'etrique sur $n$ lettres. Pour $r < n$, on identifiera $B_r$
au sous-groupe de $B_n$ form\'e des tresses qui laissent les
$n-r$ derniers brins droits. Soit $\k$ un $\Q$-anneau, et
$\mathcal{T}_n(\k)$ la $\k$-alg\`ebre de Lie d'holonomie
d\'efinie par g\'en\'erateurs $t_{ij}, 1 \leq i,j \leq n$ et relations
$t_{ii} = 0$, $t_{ij} = t_{ji}$, $[t_{ij},t_{kl}] = 0$ si $\# \{i,j,k,l\} = 4$
et $[t_{ij},t_{ik}+t_{kj}] = 0$ pour tous $i,j,k$.
Soit alors $\AN(\k)$ la version infinit\'esimale du groupe de tresses, c'est-\`a-dire
le produit semi-direct de l'alg\`ebre de groupe $\k \SN$ du groupe sym\'etrique par l'alg\`ebre
enveloppante universelle de $\mathcal{T}_n(\k)$,
l'action donnant le produit
semi-direct \'etant $s.t_{ij} = t_{s(i) s(j)}$. Ces alg\`ebres
sont naturellement gradu\'ees, par $deg(t_{ij}) = 1$ et $deg(s) = 0$ si
$s\in \SN$.
On note
$\widehat{\mathcal{T}_n}(\k)$ et $\widehat{\AN}(\k)$ leur compl\'et\'e
par rapport \`a cette graduation.

On notera $\sigma_1,\dots,\sigma_{n-1}$ les g\'en\'erateurs d'Artin
classiques de $B_n$, $\xi_{ij} = \sigma_{j-1} \dots \sigma_{i+1} \sigma_i^2 
\sigma_{i+1}^{-1} \dots \sigma_{j-1}^{-1}$ les g\'en\'erateurs 
traditionnels de $P_n$, et
$$\delta_r = \sigma_{r-1} \dots \sigma_2 \sigma_1^2
\sigma_2 \dots \sigma_{r-1}.$$
Rappelons que ces \'el\'ements $\delta_2,\dots
,\delta_r$ engendrent un sous-groupe ab\'elien libre \`a l'int\'erieur
de $P_n$, et que de plus $\delta_r = \gamma_r \gamma_{r-1}^{-1}$
o\`u $\gamma_r = (\sigma_1 \dots \sigma_{r-1})^r \in B_r$ engendre le centre
de $B_r$ pour $r \geq 3$. Les analogues infinit\'esimaux de $\delta_r$
et $\gamma_r$ sont not\'es $Y_r = t_{1r} + \dots + t_{r-1,r}$ et
$Z_r = \sum_{1 \leq i,j \leq r} t_{ij}$. On a $Y_r = Z_r - Z_{r-1}$.

Soit $C^r P_n$ pour $r \geq 0$ la suite centrale descendante de $P_n$,
d\'efinie par $C^0 P_n = P_n$, $C^1 P_n = [P_n,P_n]$ (groupe des
commutateurs) et $C^{r+1} P_n = [P_n,C^r P_n]$. On note $P_n(\k)$ ce
que Drinfeld appelle la
com\-pl\'e\-tion $\k$-pro\-uni\-po\-ten\-te de $P_n$, c'est-\`a-dire la
limite projective des $(P_n/C^r P_n) \otimes \k$. Comme la suite
centrale descendante de $P_n$ est s\'eparante, c'est-\`a-dire que les
sous-groupes $C^r P_n$ ont une intersection triviale, et que de
plus les quotients $P_n/C^r P_n$ sont sans torsion, le morphisme naturel
de $P_n$ dans ces compl\'etions est injectif. L'action par conjugaison de $B_n$
sur son sous-groupe $P_n$ s'\'etend naturellement en une action de
$B_n$ sur $P_n(\k)$. On note $B_n(\k)$ le quotient de $B_n \ltimes P_n(\k)$
par le sous-groupe (distingu\'e) engendr\'e par les \'el\'ements $\sigma. \sigma^{-1}$
pour $\sigma \in P_n = B_n \cap P_n(\k)$.

\subsection{Groupe libre et groupe de Hausdorff}

Soit $\k$ un corps de caract\'eristique 0,
$\mathcal{A}(\k) = \k \ll x,y \gg$ la $\k$-alg\`ebre des s\'eries formelles
\`a coefficients dans $\k$ en des ind\'etermin\'ees non commutatives
$x$ et $y$. Si $\k$ est un corps topologique, on munit $\mathcal{A}(\k)$
de la topologie de la convergence simple de ses coefficients, ce qui en
fait une $\k$-alg\`ebre topologique compl\`ete. On note $\mathcal{L}(\k)$ l'ensemble des
s\'eries de Lie formelles en $x$ et $y$ \`a coefficients dans $\k$,
naturellement identifi\'ee \`a la sous-alg\`ebre de Lie de $\mathcal{A}(\k)$
form\'ee des s\'eries de Lie en $x$ et $y$, et $\mathcal{A}_+(\k) \subset
\mathcal{A}(\k)$ l'ensemble
des \'el\'ements de terme constant nul. Inversement, $\mathcal{A}(\k)$ s'identifie
\`a la compl\'etion de l'alg\`ebre enveloppante universelle de l'alg\`ebre
de Lie libre sur $x,y$ pour la graduation $\deg(x) = \deg(y) =1$. L'ordre
$\omega(f)$ de $f \in \mathcal{A}(\k)$ est par d\'efinition le
degr\'e de son mon\^ome de plus bas degr\'e.
On note encore $\mathcal{M}(\k)$ ce que Bourbaki appelle le groupe
de Magnus, c'est-\`a-dire le groupe des \'el\'ements de $\mathcal{A}(\k)$
de terme constant \'egal \`a 1, et $\mathcal{H}(\k) = \exp \mathcal{L}(\k) \subset \mathcal{M}(\k)$,
c'est-\`a-dire le groupe de Hausdorff, ensemble des \'el\'ements de
$\mathcal{A}(\k)$ qui sont grouplike pour sa structure naturelle
d'alg\`ebre de Hopf compl\'et\'ee. Si $\k$ est un corps topologique,
toutes ces parties de $\mathcal{A}(\k)$ sont munies de la topologie induite.

Il est classique (cf. \cite{BOURB} ch. II \S 5) que l'exponentielle
et le logarithme fournissent des isomorphismes entre $\mathcal{A}_+(\k)$ et
$\mathcal{M}(\k)$ d'une part, $\mathcal{L}(\k)$ et $\mathcal{H}(\k)$ d'autre part.
On v\'erifie facilement que, pour $f \in \mathcal{A}_+(\k)$, les coefficients
de $\exp(f)$ sont des polyn\^omes en les coefficients de $f$,
et qu'ainsi l'application $\exp : \mathcal{A}_+(\k) \to
\mathcal{M}(\k)$ est continue si $\k$ est un corps topologique.
Il en est de m\^eme pour l'application logarithme, ce qui montre que
$\mathcal{A}_+(\k)$ et
$\mathcal{M}(\k)$ d'autre part, $\mathcal{L}(\k)$ et $\mathcal{H}(\k)$ d'autre part,
sont hom\'eomorphes par ces applications.

Si $f(x,y) \in \mathcal{M}(\k)$ et $u,v \in \mathcal{L}(\k) \subset \mathcal{A}_+(\k)$,
on d\'efinit classiquement la substitution $f(u,v)$.
Si $f(x,y) \in \mathcal{M}(\k)$ et $u,v \in \mathcal{H}(\k) = \exp \mathcal{L}(\k)$,
on d\'efinit, suivant l'usage de Drinfeld dans \cite{DRIN},
$f(u,v)$ comme $f(\log(u),\log(v))$. Cette \'ecriture s'\'etend
naturellement et permet de d\'efinir $f(u,v)$ pour $u,v$ dans la compl\'etion
$\k$-pro-unipotente de n'importe quel groupe. En particulier,
pour $u,v \in P_n(\k)$, $f(u,v) \in P_n(\k)$.

Soit $F$ le groupe libre en deux g\'en\'erateurs $x,y$. Il est classique
(cf. \cite{BOURB} ch. II \S 5 no. 3, th. 1) que l'application
$x \mapsto e^x$, $y \mapsto e^y$ se prolonge en un morphisme de
groupes qui plonge $F$ dans $\mathcal{H}(\k)$, et que ce plongement
se factorise par la compl\'etion pro-nilpotente $\widehat{F}$ de
$F$ : on a $F \hookrightarrow \widehat{F} \hookrightarrow \mathcal{H}(\k)$,
le deuxi\`eme morphisme \'etant continu si $\k$ est un corps topologique.
De m\^eme, l'application $x \mapsto 1+x$, $y \mapsto 1+y$
plonge $F$ et $\widehat{F}$ dans $\mathcal{M}(\k)$ de fa\c con continue. Si
l'on \'etend les notations $\mathcal{M}(\k)$ et $\mathcal{A}(\k)$
de fa\c con naturelle au cas o\`u $\k$ est un anneau unitaire,
cette derni\`ere application plonge en fait $F$ et $\widehat{F}$ dans
$\mathcal{M}(\Z)$.

\subsection{Associateurs de Drinfeld}

Pour tout $\Q$-anneau $\k$ et tout $\la \in \k$, on d\'efinit
suivant \cite{DRIN} l'ensemble $\Ass_{\la}(\k)$ des $\Phi \in \mathcal{A}(\k)$ qui satisfont
les \'equations suivantes :
\begin{eqnarray}
\Delta(\Phi) = \Phi \hat{\ot} \Phi \label{grouplike} \\
\Phi(y,x) = \Phi(x,y)^{-1} \label{inverse} \\
e^{\la x/2} \Phi(z,x) e^{\la z/2} \Phi(y,z) e^{\la y/2} \Phi(x,y) = 1 \label{hexagone} \\
(d_3 \Phi)(d_1 \Phi) = (d_0 \Phi)(d_2 \Phi)(d_4 \Phi) \label{pentagone} 
\end{eqnarray}
avec $z = -x-y$. Les \'equation (\ref{hexagone}) et (\ref{pentagone}) sont appel\'ees
\'equation de l'hexagone et du pentagone, respectivement.
Dans la relation (\ref{grouplike}), le symbole $\hat{\otimes}$ d\'esigne le
produit tensoriel compl\'et\'e et $\Delta$ le coproduit associ\'e \`a l'identification de $\mathcal{A}(\k)$
avec la big\`ebre enveloppante de $\mathcal{L}(\k)$.
Finalement, pour que l'\'equation (\ref{pentagone}) ait un sens, il nous faut
d\'efinir
$$
\left\lbrace \begin{array}{lcl}
d_3 \Phi & = & \Phi(t_{12},t_{23} + t_{24} ) \\
d_1 \Phi & = & \Phi(t_{13} + t_{23},t_{34}) \\
d_0 \Phi & = & \Phi(t_{23},t_{34}) \\
d_2 \Phi & = & \Phi(t_{12} + t_{13},t_{24} + t_{34}) \\
d_4 \Phi & = & \Phi(t_{12} ,t_{23}) \\
\end{array} \right.
$$
de sorte que (\ref{pentagone}) est une \'equation dans $\widehat{\mathsf{U} \mathcal{T}}_4(\k)$.
Remarquons finalement que l'\'equation (\ref{grouplike}) est \'equivalente
\`a dire que $\Phi \in \mathcal{H}(\k)$. On montre facilement que,
pour tout $\Q$-anneau $\k$ et tout $\la \in \k$, $\Ass_{\la}(\k) = 
\Ass_{-\la}(\k)$. En particulier, \`a tout $\Phi \in \Ass_{\la}(\k)$
on peut associer un autre associateur $\bar{\Phi}(x,y) = \Phi(-x,-y) 
\in \Ass_{\la}(\k)$. Si l'on suppose qu'existe un $\Phi \in \Ass_{\la}(\k)$,
il est facile de d\'eterminer la forme de ses premiers termes.
En particulier, il existe toujours $\cc \in \k$ tel que
$\Phi(x,y)$ vaille $1 + \frac{\la^2}{24} [x,y] + \cc \left( [x,[x,y]] - 
[y,[y,x]] \right)$ plus des termes d'ordre sup\'erieur (cf. par exemple
\cite{ASSOC} prop. 1). \'Etant donn\'ee l'importance
pour nous de ce coefficient, nous le consid\'erons comme une
fonction $\cc(\Phi)$ de l'associateur.

Dans \cite{DRIN}, Drinfeld construit explicitement un associateur
$\Phi_{KZ} \in \Ass_{1}(\C)$ \`a partir de l'\'etude du syst\`eme 
diff\'erentiel de Knizhnik et Zamolodchikov. 
Une formule explicite pour les coefficients
de $\Phi_{KZ}$ est due \`a Le et Murakami. En particulier, $\cc(\Phi_{KZ})
 = 
-\zeta(3)/(2\ii \pi)^3 \neq 0$~, donc
$\Phi_{KZ} \neq \bar{\Phi}_{KZ}$ puisque $\cc(\bar{\Phi}) = -\cc(\Phi)$.
Il introduit d'autre part l'ensemble des
associateurs pairs $\Ass_{\la}^0(\k)$, d\'efinis comme l'ensemble des
$\Phi \in \Ass_{\la}(\k)$ tels que $\bar{\Phi} = \Phi$,
et montre $\Ass_1^0(\Q) \neq 0$ (cf. \cite{DRIN} prop. 5.3).

Ces associateurs permettent de d\'efinir des morphismes de $B_n(\k)$
vers le groupe $(\widehat{\AN}(\k))^{\times}$ des inversibles
de $\widehat{\AN}(\k)$. Plus pr\'ecis\'ement, \`a tout $\Phi \in \Ass_{\la}(\k)$
on associe un morphisme $\widetilde{\Phi} : B_n(\k) \to 
(\widehat{\AN}(\k))^{\times}$ tel que 
$$\widetilde{\Phi}(\sigma_i) = \Phi(t_{i,i+1},Y_i)
s_i \exp(\la t_{i,i+1}/2) \Phi(Y_i,t_{i,i+1})$$

On v\'erifie ais\'ement (cf. \cite{ASSOC} p. 9 prop. 2) que
$\widetilde{\Phi}(\delta_r) = \exp (\la Y_r)$ pour tout $\Phi \in \Ass_{\la}
(\k)$ et $r \in [2,n]$.

\subsection{Le groupe de Grothendieck-Teichm\"uller $\k$-pro-unipotent}

Le groupe $GT(\k)$ est l'ensemble des couples $(\la ,f) \in \kt \times
\mathcal{H}(\k)$ qui v\'erifient $f(u,v) = f(v,u)^{-1}$
pour tous $u,v \in \mathcal{H}(\k)$,
$f(w,u)w^mf(v,w)v^m f(u,v)u^m = 1$ pour tous $u,v,w \in \mathcal{H}(\k)$
tels que $uvw=1$ et $m = (\la -1)/2$ (\'equation de l'hexagone),
et enfin une \'equation dans $P_4(\k)$ :
$$
f(\xi_{12},\xi_{23} \xi_{24}) f(\xi_{13}\xi_{23},\xi_{34}) = 
f(\xi_{23},\xi_{34}) f(\xi_{12} \xi_{13},\xi_{24}\xi_{34})
f(\xi_{12},\xi_{23}).
$$
Si $\k$ est un corps topologique, on munit $GT(\k)$ de la topologie
naturelle de $\kt \times \mathcal{H}(\k)$. C'est un groupe (topologique)
muni de la loi $(\la_1,f_1).(\la_2,f_2) = (\la,f)$ avec $\la = \la_1 \la_2$
et
$$
f(u,v) = f_1(f_2(u,v)u^{\la_2} f_2(u,v)^{-1},v^{\la_2}) f_2(u,v)
$$
En particulier, on a un morphisme de groupes (topologiques) $GT(\k) \to
\kt$, dont le noyau est not\'e $GT_1(\k)$. Ces groupes admettent des
analogues infinit\'esimaux, d\'efinis comme suit. Soit 
$GRT_1(\k) = \Ass_0(\k)$, muni de la loi $f_1 . f_2 = f$
avec
$f(x,y) = f_1(f_2(x,y)xf_2(x,y)^{-1},y)f_2(x,y)$. C'est une loi
de groupe (topologique) et on a une action de $\kt$ sur $GRT_1(\k)$,
d\'efinie par $c.f(x,y) = f(c^{-1}x,c^{-1}y)$. L'analogue infinit\'esimal
de $GT(\k)$ est alors d\'efini comme $GRT(\k) = \kt \ltimes GRT_1(\k)$. Plus
g\'en\'eralement, la m\^eme formule donne une
action \`a droite de $GRT_1(\k)$ sur $\Ass_{\la}(\k)$, pour tout
$\la \in \k$ : si $\Phi \in \Ass_{\la}(\k)$ et $f \in \Ass_0(\k) = GRT_1(\k)$,
on a $\Phi.f \in \Ass_{\la}(\k)$. De plus, cette action est libre
et transitive d'apr\`es \cite{DRIN} prop. 5.5.

Le lien entre les groupes $GRT_1(\k)$ et $GT_1(\k)$ est donn\'e par leur 
action commune sur les associateurs. On a en effet une action \`a gauche
$(f,\Phi)\mapsto f.\Phi$ de $GT_1(\k)$ sur $\Ass_{\la}(\k)$, o\`u
$$
(f.\Phi)(x,y) = f(\Phi(x,y) \exp(x) \Phi(x,y)^{-1},\exp(y)) \Phi(x,y).
$$
Comme cette action est \'egalement libre et transitive d'apr\`es \cite{DRIN} prop. 5.1
et que les deux actions commutent,
\`a tout $\Phi \in \Ass_{\la}(\k)$ est associ\'e un isomorphisme
$\iota_{\Phi} : GRT_1(\k)
\to GT_1(\k)$ d\'efini par $\Phi.f  = \iota_{\Phi}(f).\Phi$
pour tout $f \in GRT_1(\k)$.

Le groupe $GT_1(\k)$ agit \`a droite sur $B_n(\k)$ par les formules
$\sigma_1 \mapsto \sigma_1$, $\sigma_r \mapsto f(\delta_r,\sigma_r^2)^{-1}
\sigma_r f(\delta_r,\sigma_r^2)$ pour $f \in GT_1(\k)$. De plus, pour tous
$\Phi \in \Ass_1(\k)$, $\sigma \in B_n$ et $f \in GT_1(\k)$, on a
$\widetilde{\Phi}(\sigma . f) = \widetilde{(f.\Phi)}(\sigma)$. Effet,
il suffit de v\'erifier cette \'egalit\'e sur les g\'en\'erateurs d'Artin
parce que $GT_1(\k)$ agit sur $B_n(\k)$ par automorphismes de groupe ; or
$\widetilde{(f.\Phi)}(\sigma_r) = Q \widetilde{\Phi}(\sigma_r) Q^{-1}$
avec $Q$ \'egal \`a
$$
f\left( \Phi(t_{r,r+1},Y_r) \exp(t_{r,r+1}) \Phi(Y_r,t_{r,r+1}), \exp\, Y_r
\right)
= f\left( \widetilde{\Phi}(\sigma_r^2),\widetilde{\Phi}(\delta_r)\right) .$$
On en d\'eduit l'\'egalit\'e voulue.

\subsection{Un r\'esultat de densit\'e}

Le but de ce paragraphe est de d\'emontrer le th\'eor\`eme suivant.

\begin{theo} \label{theodensite} Si $L$ est un corps topologique et $L' \subset L$
est un sous-corps topologique dense de $L$, alors $GT_1(L')$ est dense
dans $GT_1(L)$.
\end{theo}

Pour d\'emontrer le th\'eor\`eme, on choisit $\Phi \in \Ass_1(\Q) \subset 
\Ass_1(L') \subset \Ass_1(L)$ et on utilise l'isomorphisme $\iota_{\Phi}$
entre $GRT_1(L)$ et $GT_1(L)$. Il envoie $GRT_1(L')$ sur
$GT_1(L')$. Ces isomorphismes respectent la topologie, comme le
montre le lemme suivant.

\begin{lemme} Soit $\k$ un corps topologique. La bijection $\iota_{\Phi}$
associ\'ee \`a $\Phi \in \Ass_1(\k)$ entre $GRT_1(\k)$ et $GT_1(\k)$
est un hom\'eomorphisme.
\end{lemme}
\begin{proof}
Soit $g \in GRT_1(\k)$ et $f = \iota_{\Phi}(g) \in GT_1(\k)$.
Par d\'efinition, la relation entre $f$ et $g$ est donn\'ee par
les formules
\begin{eqnarray}
\Phi'(x,y)  & =  & \Phi(g(x,y) x g(x,y)^{-1},y) g(x,y) \label{grtass} \\
\Phi'(x,y)  & =  & f(\Phi(x,y) \exp(x) \Phi(x,y)^{-1}, \exp(y)) \Phi(x,y)
\label{assgt}
\end{eqnarray}
Il est clair que la relation (\ref{grtass}) associe continument, \`a tout
\'el\'ement $g \in \mathcal{M}(\k)$, un \'el\'ement $\Phi' \in
\mathcal{M}(\k)$. Posant $f = \exp \circ F$ avec $F \in \mathcal{L}(\k)$,
la relation (\ref{assgt}) dit qu'alors
$$
F\left( \log \left( \Phi(x,y) \exp(x) \Phi(x,y)^{-1} \right),y \right)
 = \log \left( \Phi'(x,y) \Phi(x,y)^{-1} \right).
$$
Or il existe $P \in \mathcal{A}(\k)$, $\omega(P) \geq 2$, tel que
$\log \left( \Phi(x,y) \exp(x) \Phi(x,y)^{-1} \right)$ vaille $x + P$. Comme
l'application $\Phi' \mapsto \log \left( \Phi' \Phi^{-1} \right)$ est
continue, il suffit de montrer que l'endomorphisme continu $\Delta$
du $\k$-espace vectoriel topologique $\mathcal{A}_+(\k)$ d\'efini par 
$F(x,y) \mapsto F(x+P,y)$ admet une r\'eciproque continue. Or, pour tout
$F \in \mathcal{A}_+(\k)$, $(\Delta -\Id)(F)$ est d'ordre strictement
sup\'erieur \`a $\omega(F)$, ainsi l'ordre de $(\Delta - \Id)^n(F)$
est au moins $\omega(F)+n$ et $\Delta^{-1} = \sum (\Delta- \Id)^n (-1)^n$
est continu pour la topologie de la convergence simple puisque chaque
$(\Delta-\Id)^n$ l'est. La continuit\'e du morphisme inverse $GT_1(\k) \to
GRT_1(\k)$ se montre de m\^eme.  
\end{proof}

Il suffit donc de montrer que $GRT_1(L')$ est dense dans $GRT_1(L)$.
Or, pour tout corps $\k$, $GRT_1(\k)$ est naturellement d\'efini comme
l'ensemble des $\k$-points d'un $\Q$-sch\'ema en groupe pro-unipotent
$GRT_1$ dont l'alg\`ebre de Lie $\mathfrak{grt}_1$ est une sous-alg\`ebre de Lie
de $\mathcal{L}$ d\'efinie par
des \'equations lin\'eaires (cf. \cite{DRIN} p. 851 et prop. 5.7)
\`a coefficients rationnels. On en d\'eduit $\mathfrak{grt}_1(L) = 
\mathfrak{grt}_1(L') \otimes_{L'} L$, et $\mathfrak{grt}_1(L')$ est dense dans $\mathfrak{grt}_1(L)$.
Pour avoir la conclusion du th\'eor\`eme, il suffit de montrer
que l'application exponentielle de cette alg\`ebre de Lie
dans son groupe de Lie est continue. Pour tout $f \in \mathcal{A}_+(\k)$,
notons $D_f$ l'unique d\'erivation continue de $\mathcal{A}(\k)$
qui v\'erifie $D_f(x) = [f,x]$ et $D_f(y) = 0$, et $s_{f}$ l'endomorphisme
lin\'eaire de $\mathcal{A}(\k)$ qui \`a $g \in \mathcal{A}(\k)$ associe
$s_f(g) = gf + D_f(g)$. On a $\omega(s_f(g)) > \omega(g)$. Comme $GRT_1$
est pro-unipotent en tant que sch\'ema en groupe sur $\Q$, l'application
exponentielle associ\'ee est surjective. L'image de $f \in \mathfrak{grt}_1(\k)$
dans $GRT_1(\k)$ est donn\'ee par $\exp(s_f)(1)$ (cf. \cite{DRIN}, ou
\cite{JOE} prop. 2.9 pour plus de d\'etails), on en d\'eduit que cette
application est continue. Le th\'eor\`eme est ainsi d\'emontr\'e.

\subsection{Caract\`eres de Soul\'e}
Pour $l$ un nombre premier, notons $\mu_{l^{\infty}}$ la limite inductive
sur $n$ des $\mu_{l^n}$ et $\Gamma$ le groupe de Galois
absolu de $\Q(\mu_{l^{\infty}})$, muni de sa topologie naturelle de groupe
profini.

Soul\'e a d\'efini une classe de caract\`eres de $\Gamma$ \`a valeurs
dans le groupe additif $\Z_l$, dont nous rappelons bri\`evement la
construction. Soit $m \geq 3$ impair. Pour tout $n \geq 1$ (resp. $n \geq 2$ si $l = 2$)
on note $\zeta_n$ une racine primitive $l^n$-i\`eme de 1, et
$$
\eps_{m,n} = \prod_a (\zeta_n^a-1)^{[a^{m-1}]}
$$
o\`u $a$ parcourt les entiers strictement compris entre 0 et $l^n$
qui sont premiers \`a $l$ et $[a^{m-1}]$ d\'esigne le reste
de la division euclidienne de $a^{m-1}$ par $l^n$. Comme $m-1$ est pair,
$\eps_{m,n}$ est totalement r\'eel et totalement positif. Pour $u,v >0$ on note $u^v = \exp(v \log u)$.
Il existe alors, pour tout $\sigma \in \Gamma$,
un unique $\kappa_m(\sigma) \in \Z_l$ tel que, pour tout $n \geq 1$
($n \geq 2$ si $l=2$),
$$
\sigma( (\eps_{m,n}^{\frac{1}{l^n}}) = \sigma(\eps_{m,n})^{\frac{1}{l^n}}
\zeta_n^{\kappa_m(\sigma)}
$$
Ces applications $\kappa_m$ sont les caract\`eres de Soul\'e.

Si $l=2$, $\eps_{3,2} = \eps_{3,3} = \eps_{5,3} = \eps_{5,2} = 2$,
donc $\kappa_3(\sigma)$ et $\kappa_5(\sigma)$ sont congrus modulo 8.
D'autre part, le fait que $\sqrt{2} \in \Q(\mu_8)$ implique
que $\kappa_3(\sigma)$ est divisible par 2.
Si $l = 3$ on a $\eps_{3,1} = \eps_{5,1} = (j-1)(j^2-1) = 3$,
donc de m\^eme $\kappa_3(\sigma)$ est congru \`a $\kappa_5(\sigma)$ modulo 3.

On utilisera \'egalement la notation, \`a $l$ fix\'e, $\kappa_m^*
(\sigma) = \kappa_m(\sigma)/(l^{m-1}-1).$ Le principal r\'esultat
dont nous aurons besoin concernant ces caract\`eres est le suivant :
\begin{prop} \label{k3nonnul}
Pour tout nombre premier $l$, $\kappa_3 \neq 0$.
\end{prop}
\begin{proof} Si $l$ est impair, cela d\'ecoule du
r\'esultat g\'en\'eral de Soul\'e sur la non trivialit\'e
des $\kappa_m$ (cf. \cite{SOULE1,SOULE2}, voir \'egalement \cite{ICHI}).
Le cas $l = 2$ est \'el\'ementaire : il suffit en effet de
montrer que $X^4-2$ n'admet pas de racine dans $L = \Q(\mu_{2^{\infty}})$,
extension galoisienne ab\'elienne de $\Q$. C'est le cas, soit
parce que $\Q(\sqrt[4]{2}) \subset \R$ serait autrement
une extension galoisienne de $\Q$ (car Gal($L | \Q$) aurait Gal($\Q(
\sqrt[4]{2})|\Q$) comme sous-groupe, n\'ecessairement distingu\'e),
soit parce que $\Q(\sqrt[4]{2},\ii)$, extension non ab\'elienne de
$\Q$, serait autrement incluse dans $L$.
\end{proof}

On a un morphisme $\Gal(\bar{\Q}|\Q) \to GT(\Q_l)$ (cf. \cite{DRIN,IHARA}).
Son compos\'e avec le morphisme naturel $GT(\Q_l) \to \Q_l^{\times}$
est la composante en $l$ du caract\`ere cyclotomique. On en
d\'eduit un morphisme $m$ du groupe de Galois absolu
de $\Q(\mu_{l^{\infty}})$ vers $GT_1(\Q_l)$. Pour all\'eger les notations,
pour $\sigma \in \Gal( \bar{\Q}|\Q(\mu_{l^{\infty}}))$
nous noterons $\sigma.\Phi$, $\chi(\sigma)$ pour $m(\sigma).\Phi$,
$\chi \circ m (\sigma)$.

\subsection{$GT_1(\k)$ et $GRT_1(\k)$ en degr\'e 3}

On note $w_x = [x,[x,y]]$ et $w_y = [y,[y,x]]$. Soient $\la \in \kt$,
$\Phi \in \Ass_{\la}(\k)$, $f \in GRT_1(\k)$, $F \in GT_1(\k)$. Pour
$a,b \in \mathcal{A}(\k)$, on note $a \equiv b$ si $\omega(a-b) \geq 4$.
On sait qu'existe $\cc \in \k$ tel que $\Phi \equiv 1 + \frac{\la^2}{24}
[x,y] + \cc(w_x - w_y)$. Comme $GRT_1(\k) = \Ass_0(\k)$, il existe
$z \in \k$ tel que $f \equiv 1 + z(w_x - w_y)$. Comme
$F \in \mathcal{H}(\k)$, $F = \exp \psi$ pour un certain $\psi \in
\mathcal{L}(\k)$. On a donc $\psi \equiv a_1 x + a_2 y + a_3 [x,y]
+ a_4 w_x + a_5 w_y$ pour certaines valeurs des $a_i \in \k$. L'\'equation
$F(u,v) F(v,u) = 1$ implique
$a_2 = -a_1$ en degr\'e 1, puis $a_1 = 0$ en degr\'e 2, l'\'equation
de l'hexagone en degr\'e 2 implique $a_3 = 0$, et enfin la
premi\`ere \'equation en degr\'e 3 implique $a_5 + a_4 = 0$ : il
existe donc $z' \in \k$ tel que $F \equiv 1+z'(w_x - w_y)$. Les valeurs
$z , \la, \cc,z'$ \'etant fix\'ees,
$$
\begin{array}{lclcl}
(\Phi.f) (x,y) & = & \Phi(fxf^{-1},y)f & \equiv & 1 + \frac{\la^2}{24}[x,y]
+ (\cc +z)(w_x - w_y) \\
(F.\Phi)(x,y) & = & F( \Phi e^x \Phi^{-1},e^y) \Phi & \equiv & 1 +
\frac{\la^2}{24} [x,y] + (\cc + z')(w_x - w_y) \\
\end{array}
$$
L'isomorphisme $\iota_{\Phi}$ associ\'e \`a $\Phi$ induit donc
toujours l'identit\'e en degr\'e 3 (c'est \'egalement une cons\'equence
\'el\'ementaire du th\'eor\`eme de Drinfeld selon lequel
le gradu\'e de l'alg\`ebre de Lie de $GT_1(\k)$ est isomorphe \`a
$\mathfrak{grt}_1(\k)$). D'autre part, $\mathfrak{grt}_1(\k)$ est une $\k$-alg\`ebre de Lie
gradu\'ee qui contient l'\'el\'ement $w_x - w_y$ donc, pour tout
$z \in \k$, $GRT_1(\k)$ contient l'exponentielle (au sens du groupe)
de $z (w_x - w_y)$, qui est congrue \`a $ 1 + z (w_x - w_y)$. En
particulier, on d\'eduit des formules pr\'ec\'edentes
et de l'existence, pour tout $\la \in \k$, d'un associateur
$\Phi \in \Ass_{\la}(\k)$ le lemme suivant.

\begin{lemme} \label{ccnonnul} Pour tous $\cc,\la \in \k$, il existe
$\Phi \in \Ass_{\la}(\k)$ tel que $\cc(\Phi) = \cc$. En particulier,
pour tout $\la \in \k$ il existe $\Phi \in \Ass_{\la}(\k)$
tel que $\cc(\Phi) \neq 0$.
\end{lemme}

Nous allons \'egalement d\'eterminer le terme de degr\'e 3 pour
l'image du groupe de Galois. 
On a un morphisme du groupe de Galois absolu de
$\Q(\mu_{l^{\infty}})$ vers le groupe de Grothendieck-Teichm\"uller
pro-$l$, qui envoie un \'el\'ement $\sigma$ du groupe de Galois sur
un \'el\'ement $f^{\sigma}$ du compl\'et\'e pro-l $F_l$ du groupe libre
en deux g\'en\'erateurs $u,v$. Il est classique que $F_l$
se plonge dans $\mathcal{H}(\Q_l)$ par l'application $u \mapsto e^x$,
$v \mapsto e^y$. L'image de $f^{\sigma}$ par ce morphisme sera
not\'ee $Df^{\sigma}$. On a $Df^{\sigma} \in GT_1(\Q_l)$. 

On peut \'egalement plonger $F_l$ dans $\mathcal{A}(\Z_l)$
en envoyant $u$ et $v$ respectivement sur $1+x$ et $1+y$. On note
$If^{\sigma} \in \mathcal{A}(\Z_l)$ l'image de $f^{\sigma}$ par
ce morphisme. On a imm\'ediatement, dans $\mathcal{A}(\Q_l)$,
$$
Df^{\sigma}(\log(1+x), \log(1+y)) = I f^{\sigma}(x,y), \ \ 
I f^{\sigma}(e^x-1,e^y-1) = Df^{\sigma}(x,y).
$$

Introduisons alors, suivant Y. Ihara (cf. \cite{IHARA}), la d\'ecomposition naturelle
$\mathcal{A}(\k) = \k \oplus  \mathcal{A}(\k) x \oplus \mathcal{A}(\k) y$,
notons $p_x$ la projection sur les deux premi\`eres composante, et
$\pi_{ab}$ la restriction du morphisme d'ab\'elianisation $\mathcal{A}(\k)
\to \k[x,y]$ \`a $\k \oplus \mathcal{A}(\k)x$. Ihara montre que
$$
\psi^{\sigma}_{ab} = \pi_{ab} \circ p_x(If^{\sigma})
= \exp \left( \sum_{\stackrel{m \geq 3}{m\mbox{ impair }}}
\frac{\kappa_m^*(\sigma)}{m!} \left( (X+Y)^m - X^m - Y^m \right) \right)
$$
avec $X = \log(1+x)$, $Y = \log(1+y)$. En particulier,
$\psi^{\sigma}_{ab} \equiv 1 + \kappa_3^*(\sigma)(yx^2 + y^2x)/2$.
Alors, si $Df^{\sigma} \equiv 1 + \alpha( w_x - w_y)$,
on a $If^{\sigma} \equiv Df^{\sigma}$ et, de  $\pi_{ab} \circ p_x(w_x)
= -x^2y$ et $\pi_{ab} \circ p_x(w_y) = y^2x$,
on d\'eduit $\psi^{\sigma}_{ab} = 1 - \alpha(y^2x + yx^2)$
et $\alpha = -\kappa_3^*(\sigma)/2$. On a donc montr\'e
\begin{lemme} \label{ccgal} L'image de $\sigma \in \Gal(\bar{\Q} |
\Q(\mu_{l^{\infty}}))$ dans $GT_1(\Q_l)$ vaut
$$1 - \frac{\kappa_3^*(\sigma)}{2}( [x,[x,y]]-[y,[y,x]])
$$
plus des termes d'ordre sup\'erieurs. Pour tout
$\Phi \in \Ass_{\la}(\Q_l)$,
$$ \cc( \sigma.\Phi) = \cc(\Phi) - \frac{\kappa_3^*(\sigma)}{2}.
$$
\end{lemme}

\section{Action de $GT_1(\k)$ sur les repr\'esentations de $B_n$}

\subsection{G\'en\'eralit\'es}

Soit $\k$ un corps de caract\'eristique 0, $A = \k[[h]]$, $K= \k((h))$. 
Soit
$$
V_N(\k) = \{ R \in \Hom(B_n,GL_N(A)) \ \mid \ R(C^r P_n) \subset 1 + h^{r+1}M_N(A) \}.
$$
On a $V_N(\k) \subset \Hom(B_n,GL_N(A)) \subset \Hom(B_n,GL_N(K))$. Tout $R
\in V_N(\k)$ se prolonge naturellement en $\tilde{R} \in \Hom(B_n(\k),
GL_N(A))$. On en d\'eduit une action \`a gauche de $GT_1(\k)$ sur $V_N(\k)$ :
\`a $R \in V_N(\k)$ et $f \in GT_1(\k)$ est associ\'e
$f.R = S$ d\'efinie par $S(\sigma) = \tilde{R}(\sigma.f)$ pour tous
$\sigma \in B_n$. Soit maintenant $\mathcal{V}(\k) = \Hom(\AN(\k),M_N(\k))$
et $\Phi \in \Ass_{\la}(\k)$. A tout $\rho \in \mathcal{V}(\k)$
on associe $\bar{\rho} \in \Hom(\widehat{\AN}(\k),M_N(A))$
d\'efinie par $\bar{\rho}(t_{ij}) = h t_{ij}$,
$\bar{\rho}(s) = \rho(s)$ pour $ s \in \SN$, et $\widehat{\Phi}(\rho)
= \bar{\rho} \circ \widetilde{\Phi} \in \Hom(B_n,GL_N(A))$. On v\'erifie
ais\'ement $\widehat{\Phi}(\rho) \in V_N(\k)$. Les propri\'et\'es de
th\'eorie des repr\'esentations du foncteur $\widehat{\Phi}$
ont \'et\'e \'etudi\'ees par l'auteur dans \cite{ASSOC}.

Soit maintenant $f \in GT_1(\k)$. On a, pour tout $\sigma \in B_n$, 
$$
\begin{array}{lclcl}
\left[ f.\widehat{\Phi}(\rho) \right] (\sigma) & = &
\widehat{\Phi}(\rho)(\sigma.f) & = & \bar{\rho} \circ \widetilde{\Phi}(\sigma.f) \\
& = & \bar{\rho} \circ \widetilde{(f.\Phi)} (\sigma) & = & \widehat{(f.\Phi)}
(\rho)(\sigma) \\
& = & \widehat{\Phi.\iota_{\Phi}^{-1}(f)}(\rho)(\sigma) \\
\end{array}
$$
d'o\`u $f.\widehat{\Phi}(\rho) = \widehat{f.\Phi}(\rho) = 
\widehat{\Phi.\iota^{-1}_{\Phi}(f)}(\rho)$. On s'int\'eresse ici \`a une classe
particuli\`ere de repr\'esentations de $B_n$. Nous dirons d'une repr\'esentation $R : B_n \to GL_N(A)$ qu'elle est (absolument)
irr\'eductible s'il en est ainsi dans $GL_N(K)$. Deux telles repr\'esentations $R$ et $R'$ seront dites isomorphes si elles le sont dans
$GL_N(K)$.

\begin{defi} 
Une repr\'esentation absolument irr\'eductible $R \in V_N(\k)$ est dite \emph{GT-rigide} si $g.R$ est isomorphe \`a $R$ pour tout $g \in GT_1(\k)$.
\end{defi}

Comme l'action de $GT_1(\k)$ sur $V_N(\k)$ commute \`a l'action par conjugaison
\`a droite de $PGL_N(K)$ sur $\Hom(B_n,GL_N(K))$,
on associe \`a tout $R \in V_N(\k)$
GT-rigide une repr\'esentation projective $Q_R : GT_1(\k) \to PGL_N(K)$
par $g.R = Q_R(g)^{-1} R Q_R(g) = R.Q_R(g)$.

Soit $\rho : \AN(\k) \to M_N(\k)$ une repr\'esentation absolument
irr\'eductible de $\AN(\k)$. D'apr\`es \cite{ASSOC} la repr\'esentation
$\widehat{\Phi}(\rho)$ associ\'ee est absolument irr\'eductible pour
tout $\Phi \in \Ass_{\la}(\k)$, pourvu que $\la \in \kt$. Si $\widehat{\Phi}(\rho)$ est GT-rigide pour un tel $\Phi$, il en sera
de m\^eme pour tous. Par abus, nous dirons alors que $\rho$ est GT-rigide.

D'autre part, si $R = \widehat{\Phi}(\rho)$ est GT-rigide, cette
repr\'esentation
s'\'etend en une repr\'esentation de $GT_1(L)$ dans $PGL_N(L((h)))$
pour tout surcorps $L$ de $\k$. Enfin, pour toutes repr\'esentations
$R_1 \in V_{N_1}(\k)$, $R_2 \in V_{N_2}(\k)$, $g \in GT_1(\k)$ et $\sigma
\in B_n$, on a $R_1 \otimes R_2 \in V_{N_1N_2}(\k)$ et
$(g.(R_1 \otimes R_2))(\sigma) = \widetilde{R_1 \otimes R_2}(\sigma.g)
= \widetilde{R_1}(\sigma.g) \otimes \widetilde{R_2}(\sigma.g)$
donc $g.(R_1 \otimes R_2) = (g.R_1) \otimes (g.R_2)$. Si $R_1$ et $R_2$
sont GT-rigides, on en d\'eduit que $R_1 \otimes R_2$ est conjugu\'ee
\`a $g.(R_1 \otimes R_2)$ par $Q_{R_1} \otimes Q_{R_2}$.

\subsection{Agr\'egeance et caract\`eres}

\subsubsection{Action de $\GM(\k)$} On identifie $\GM(\k) = \k^{\times}$ \`a un sous-groupe des $\k$-automorphismes
de $K = \k((h))$, en faisant agir $\alpha \in \GM(\k)$ sur $f(h) \in
K$ par $f \mapsto f^{\alpha}$ avec $f^{\alpha}(h) = f(\alpha h)$.
L'anneau $A \subset K$ est laiss\'e stable. Si $\k$ est un corps topologique,
cette action est continue. Elle s'\'etend naturellement \`a $GL_N(K)$
et $PGL_N(K)$, donc \'egalement \`a $\Hom(G,GL_N(K))$
et $\Hom(G,PGL_N(K))$ pour tout groupe $G$. Soit $R \in V_N(\k)$. Il existe
a priori deux fa\c cons de faire agir $\alpha \in \GM(\k)$ sur $\tilde{R}$,
soit par $\widetilde{R^{\alpha}}$, soit par $(\widetilde{R})^{\alpha}$,
c'est-\`a-dire par l'action de $\GM(\k)$ soit sur $\Hom(B_n,GL_N(A))$,
soit sur $\Hom(B_n(\k), GL_N(K))$. Comme l'action de $\GM(\k)$ sur $K$ est
continue pour la topologie $h$-adique (ce qui revient \`a munir $\k$ de
la topologie discr\`ete), on v\'erifie facilement que
$\widetilde{R^{\alpha}}=(\widetilde{R})^{\alpha}$, donc que l'action
de $\GM(\k)$ est bien d\'efinie. De plus, $V_N(\k)$ est stable sous
l'action de $\GM(\k)$. On en d\'eduit que l'action de $GT_1(\k)$
sur $V_N(\k)$ commute \`a celle de $\GM(\k)$ : pour tous $g \in
GT_1(\k)$, $\alpha \in \GM(\k)$, $g.R^{\alpha} = (g.R)^{\alpha}$. 
En particulier, si $R$ est GT-rigide, $R^{\alpha}$ l'est \'egalement
pour tout $\alpha \in \GM(\k)$ et $Q_{R^{\alpha}} = (Q_R)^{\alpha}$.
D'autre part, $\GM(\k)$ agit par automorphismes sur $\AN(\k)$ : \`a $\alpha
\in \GM(\k)$ on associe l'unique automorphisme d'alg\`ebre de Hopf
tel que $t_{ij} \mapsto \alpha t_{ij}$ pour $1 \leq i,j \leq n$ et qui
laisse $\SN \subset \AN(\k)$ invariant. On en d\'eduit une action
$\rho \mapsto \rho^{\alpha}$ sur $\mathcal{V}_N(\k)$ ; elle
v\'erifie $\widehat{\Phi}(\rho^{\alpha}) = \widehat{\Phi}(\rho)^{\alpha}$
pour tout $\Phi \in \Ass_{\la}(\k)$.

\subsubsection{Repr\'esentations agr\'egeantes}

Soit $\mathcal{D}_n$ la sous-alg\`ebre de Lie commutative de
$\mathcal{T}_n$ engendr\'ee par $Y_2,\dots,Y_n$. Suivant
\cite{ASSOC} on appelle agr\'egeante toute
$\rho : \AN(\k) \to M_N(\k)$ telle que $\rho(\mathsf{U} \mathcal{D}_n)$
\'egale la sous-alg\`ebre de $M_N(\k)$ form\'ee des matrices diagonales.
Pour une re\-pr\'e\-sen\-ta\-tion a\-gr\'e\-gean\-te, l'ir\-r\'e\-duc\-ti\-bi\-li\-t\'e
\'equivaut \`a l'absolue ir\-r\'e\-duc\-ti\-bi\-li\-t\'e (cf. \cite{IRRED}, cor. 2 de la
prop. 3). De plus, si $\rho_1 \in
\mathcal{V}_{N_1}(\k)$ et $\rho_2 \in \mathcal{V}_{N_2}(\k)$
sont deux repr\'esentations agr\'egeantes et (absolument) irr\'eductibles,
il en est de m\^eme de $\rho_1^{\alpha_1} \otimes \rho_2^{\alpha_2}$
pour $(\alpha_1,\alpha_2) \in \GM(\k)^2$ g\'en\'erique d'apr\`es 
\cite{IRRED}.

Soit $D_n$ le sous-groupe ab\'elien libre de $P_n$
engendr\'e par $\delta_2,\dots,\delta_n$. Supposons que
$R \in V_N(\k)$ est GT-rigide et telle que $R(\delta_2),...,R(\delta_n)$ engendre 
la sous-alg\`ebre de $M_N(K)$ form\'ee des matrices diagonales.
Cette derni\`ere condition est en particulier satisfaite si $R = \widehat{\Phi}(\rho)$ pour
un certain $\Phi \in \Ass_{\la}(\k)$, $\la \in \kt$ et
$\rho \in \mathcal{V}_N(\k)$ agr\'egeante, puisqu'alors
$R(\delta_r) = \exp ( \la h Y_r)$ pour $r \in [2,n]$. 
D'autre part, pour tout $r \in [2,n]$ et pour tout $g \in GT_1(\k)$,
$\delta_r.g = \delta_r$. Cette \'egalit\'e peut se d\'emontrer directement
\`a partir des \'equations de d\'efinition de $GT_1(\k)$, ou
se d\'eduire de l'existence d'un $\Phi \in \Ass_1(\k)$ : on a
$\widetilde{\Phi}(\delta_r.g) = \widetilde{(f.\Phi)}(\delta_r)
= \widetilde{\Phi}(\delta_r)$, donc $\delta_r.g = \delta_r$. On en d\'eduit
$(g.R)(\delta_r) = R(\delta_r)$, donc tout $Q_R(g)$ commute \`a $R(\delta_r
)$ pour tout $r \in [2,n]$, et appartient ainsi
\`a l'image
dans $PGL_N(K)$ des matrices diagonales inversibles de $M_N(K)$.

En particulier, l'image de $GT_1(\k)$ dans $PGL_N(K)$ est commutative
donc, si $R = \widehat{\Phi}(\rho)$, $Q_R$ ne d\'epend
pas du choix de $\Phi \in \Ass_1(\k)$. En effet, si l'on pose
$\Phi' = \tau. \Phi$ pour un certain $\tau \in GT_1(\k)$
et $R' = \widehat{\Phi'}(\rho)$, pour tout $g \in GT_1(\k)$ on a alors
$Q_{R'}(g) = Q_R(\tau^{-1} g \tau) = Q_R(g)$.

\subsubsection{Propri\'et\'es des caract\`eres}

Le morphisme de $(K^{\times})^N$ vers $\GM(K)^{N-1}$ d\'efini par
$(x_1,\dots,x_N) \mapsto (x_2/x_1,x_3/x_2,\dots,x_N/x_{N-1})$ de noyau
$K^{\times}$ induit un isomorphisme $(K^{\times})^N/K^{\times} \to
\GM(K)^{N-1}$. On en d\'eduit $N-1$ caract\`eres
$GT_1(\k)
\to \GM(K)$ de $GT_1(\k)$ que l'on num\'erote $Q_{R,2},\dots,Q_{R,N}$
suivant l'ordre des vecteurs de base de $K^N$. Naturellement et
par convention, on pose $Q_{R,1} = 1$. Pour tout $\alpha \in \GM(\k)$
on a $Q_{R^{\alpha},i} = (Q_{R,i})^{\alpha}$. Si ces $N-1$ caract\`eres
sont distincts et non triviaux, nous dirons que $R$ est
\emph{sans r\'esonances}. Cette derni\`ere propri\'et\'e ne
d\'epend pas de l'ordre des vecteurs de base que l'on a choisi. 

De plus, comme $R$ est absolument irr\'eductible, il existe $\sigma \in
\k B_n$ tel que $R(\sigma) \in M_N(A)$ ait tous ses coefficients
non nuls. Pour tout $g \in GT_1(\k)$, on a alors $Q_R(g) (g.R)(\sigma)
= R(\sigma) Q_R(g)$. Si l'on note sous forme matricielle
$R(\sigma) = (x_{ij})$, $(g.R)(\sigma) = (y_{ij}(g))$ pour $1 \leq i,j \leq N$,
$Q_R(g) = \mathrm{diag}(a_1,a_2,\dots,a_n)$ avec $a_i = Q_{R,i}(g)$, cette
\'equation est \'equivalente \`a $x_{ij} a_j = a_i y_{i,j}(g)$
pour tous $i,j$. En particulier, $Q_{R,j}(g) = y_{1,j}(g)/x_{1,j}$
pour tout $j$. Si $\k$ est un corps topologique on en d\'eduit
que chaque $Q_{R,i}$ est continu, chaque $y_{1,j}$ \'etant une fonction
continue de $g \in GT_1(\k)$.

Pour $R \in V_N(\k)$, notons $\overline{R} \in \Hom(B_n,GL_N(\k))$
la r\'eduction de $R$ modulo $h$. Comme $R(P_n) \subset 1 + h M_N(A)$,
$\bar{R}$ se factorise par $\SN$. C'est la repr\'esentation du groupe
sym\'etrique associ\'ee \`a $R$. Si $R = \widehat{\Phi}(\rho)$,
$\overline{R}$ s'identifie \`a la restriction de $\rho$ \`a $\SN$. Pour la
repr\'esentation $R$ choisie ici, si $\bar{R}$ est (absolument)
irr\'eductible on peut choisir $\sigma \in \k B_n$ tel que $\bar{R}
(\sigma)$ ait tous ses coefficients non nuls, donc chaque $x_{1j}$ est
inversible dans $A$. On en
d\'eduit que chaque $Q_{R,i}$ est alors \`a valeur dans $\GM(A)$.
Nous avons montr\'e dans \cite{KZ} que cette situation
(i.e. que $\bar{R}$ est absolument irr\'eductible) a lieu
essentiellement lorsque $R$ se factorise par l'alg\`ebre
d'Iwahori-Hecke (cf. \cite{KZ}).

De fa\c con g\'en\'erale, si l'on suppose $R = \widehat{\Phi}(\rho)$
pour un certain $\Phi \in \Ass_1(\k)$ et $\rho$ absolument
irr\'eductible et agr\'egeante, les caract\`eres sont toujours \`a valeurs
dans $A_1 = \exp A_0$, o\`u $A_0 = hA$ et $A_1 = 1+hA$. En effet,
il existe alors (voir \cite{ASSOC} section 3.1.2) un \'el\'ement
$\sigma \in K B_n$ ind\'ependant du choix de $\Phi \in \Ass_1(\k)$
tel que $\widehat{\Phi}(\rho)(\sigma) = m_0 + h m$ avec $m \in M_N(A)$
et $m_0$ l'\'el\'ement de $M_N(\k)$ dont tous les coefficients
sont \'egaux \`a 1. On en d\'eduit $Q_{R,j}(g) \in A_1$ pour
tout $g \in GT_1(\k)$.

Sous certaines conditions (d' \og unitarit\'e \fg) on peut alors
montrer que le logarithme de ces caract\`eres est une s\'erie
formelle \emph{impaire} en $h$ sans terme lin\'eaire (voir en 4.5).

\subsection{Continuit\'e.}

Nous d\'emontrons que les repr\'esentations
projectives $Q_R : GT_1(\k) \to K^{\times}$ sont continues si $\k$ est un
corps topologique. Cela d\'ecoule de consid\'erations
g\'en\'erales que nous rappelons ici par manque de r\'ef\'erences
ad\'equates.

Si $K$ est un corps topologique et $\mathcal{A}$ une $K$-alg\`ebre
\`a unit\'e de type fini (munie de la topologie induite),
l'ensemble de morphismes d'alg\`ebres \`a
unit\'e $\Hom(\mathcal{A},M_N(K))$ est naturellement muni d'une topologie.
Elle est telle que si $g_1,\dots,g_r$ forment un syst\`eme quelconque
de g\'en\'erateurs de $\mathcal{A}$, l'application $R \mapsto (R(g_i))_{
i=1..r} \in M_N(K)^r$ est un hom\'eomorphisme de $\mathcal{A}$ sur son
image. Pour $R \in \Hom(\mathcal{A},M_N(K))$ on note $\mathcal{O}(R)$
son image sous l'action par conjugaison de $PGL_N(K)$. On a alors la

\begin{prop} \label{homeo} Si $R \in \Hom(\mathcal{A},M_N(K))$ est absolument
irr\'eductible, alors $P \mapsto P.R$ est un hom\'eomorphisme
de $PGL_N(K)$ sur $\mathcal{O}(R)$.
\end{prop}
\begin{proof} La bijectivit\'e d\'ecoule du lemme de Schur. Il s'agit
de montrer que la r\'eciproque est continue. D'apr\`es le th\'eor\`eme
de Burnside $R(\mathcal{A}) = M_N(K)$ donc quitte \`a quotienter
par $\Ker \, R$ on peut supposer $R$ bijective. On note $(E_{ij})$
pour $1 \leq i ,j \leq N$ la famille des matrices \'el\'ementaires
de $M_N(K)$ et $g_{ij} = R^{-1}(E_{ij})$. On plonge $\Hom(\mathcal{A},
M_N(K))$ dans $M^*_N(K)^{N^2}$ (de fa\c con $PGL_N(K)$-\'equivariante)
suivant ces g\'en\'erateurs, o\`u l'on note $M^*_N(K) = M_N(K)\setminus
\{ 0 \}$. Soit $(A^{u,v})$ pour $1 \leq u,v \leq N$ un \'el\'ement
de $M_N^*(K)^{N^2}$ avec $A^{u,v}$ matrice de terme g\'en\'eral
$(a^{u,v}_{i,j})$. Soit $P \in GL_N(K)$, $Q = P^{-1}$ de
terme g\'en\'eral $(p_{i,j})$, $(q_{i,j})$. Si, pour tous $u,v$, on
a $P E_{u,v} Q = A^{u,v}$, c'est-\`a-dire $p_{iu} q_{vj}
= a_{ij}^{uv}$ pour tous $i,j,u,v$, il existe $i_0,j_0,u_0,v_0$
tels que $\beta = a_{i_0 j_0}^{u_0 v_0} \neq 0$, d'o\`u
$p_{i_0 u_0} \neq 0$, $q_{v_0 j_0} \neq 0$. La restriction de l'inverse
de $P \mapsto P.R$ \`a l'ouvert de $\mathcal{O}(R)$ d\'efini par
$a_{i_0 j_0}^{u_0 v_0} \neq 0$ est alors donn\'e par la classe
dans $PGL_N(K)$ de $(a_{i j_0}^{u v_0})_{i,u}$ puisqu'alors
$a_{i,j_0}^{u v_0} / \beta = p_{i u} / p_{i_0 u_0}$, et est donc continue.
\end{proof}
On en d\'eduit imm\'ediatement :
\begin{ccor} Si $R$ est une repr\'esentation GT-rigide de $B_n$,
$Q_R : GT_1(\k) \to PGL_N(K)$ est continue.
\end{ccor}

\subsection{Trivialit\'e.}
Nous donnons une caract\'erisation des
re\-pr\'e\-sen\-ta\-tions GT-rigides $R$ pour lesquelles $Q_R$ est
triviale.

\begin{prop} \label{trivia} Soit $\rho$ une repr\'esentation GT-rigide de
$\mathfrak{B}_n(\k)$, $n \geq 3$, $\Phi \in \Ass_{\la}(\k)$, $\la \neq 0$,
$R = \widehat{\Phi}(\rho)$. Les conditions suivantes sont
\'equivalentes
\begin{itemize}
\item[(i)] $Q_R$ est triviale.
\item[(ii)] $\rho([t_{12},t_{23}]) = 0$.
\item[(iii)] $R([P_n,P_n])= \{ 1 \}$.
\end{itemize}
Si $\k = \Q_l$, ces conditions sont encore \'equivalentes \`a \\
\indent (iv) $\forall \sigma \in \Gal(\bar{\Q}|\Q(\mu_{l^{\infty}})) \ \ 
Q_R(\sigma) = 1$. 
\end{prop}
\begin{proof}
L'\'equivalence $(ii) \Leftrightarrow (iii)$ est d\'emontr\'ee
dans \cite{ASSOC}, lemme 5. $(ii) \Rightarrow (i)$ d\'ecoule
du fait qu'alors $\widehat{\Phi}'(\rho) = \widehat{\Phi}(\rho)$
pour tout $\Phi' \in \Ass_{\la}(\k)$, parce que l'image par
$\rho$ de $\mathcal{T}_n(\k)$ est alors commutative. 
On montre $(i)
\Rightarrow (ii)$ \`a partir du calcul explicite (cf. \cite{ASSOC} prop. 3)
$$
\widehat{\Phi}_2(\rho)(\sigma_2) - \widehat{\Phi}_1(\rho)(\sigma_2)
 = \rho(s_2) (\cc(\Phi_1) - \cc(\Phi_2)) \rho([t_{23},[t_{23},t_{12}]])
$$
plus des termes d'ordres sup\'erieurs, pour tous $\Phi_1,\Phi_2 \in
\Ass_{\la}(\k)$.
D'apr\`es le lemme \ref{ccnonnul} on
peut en effet supposer $\cc(\Phi) \neq 0$. Alors $\bar{\Phi} \in \Ass_{\la}(k)$
et $\cc(\bar{\Phi}) = -\cc(\Phi) \neq \cc(\Phi)$. Si $Q_R$
est triviale, $\widehat{\Phi}(\rho)(\sigma_2) = \widehat{\bar{\Phi}}(
\rho)(\sigma_2)$. On a alors
$\rho([t_{23},[t_{23},t_{12}]]) = 0$. Comme
la restriction de $\rho$ \`a $\mathcal{T}_3 = \mathcal{T}_3(\k)$ est
$\mathfrak{S}_3$-\'equivariante on a aussi
$\rho([t_{12},[t_{12},t_{23}]]) = 0$. Soit $\mathcal{T}'$ la sous-alg\`ebre
de Lie de $\mathcal{T}_3$ engendr\'ee par $t_{12}$ et $t_{23}$. On a
$\mathcal{T}_3 = \mathcal{T}' \oplus \k T$ avec $T = t_{12}+t_{13} + t_{23}$,
$\k T = Z(\mathcal{T}_3)$, $\mathcal{T}' \simeq \mathcal{T}_3/
Z(\mathcal{T}_3)$. La restriction de $\rho$ \`a $\mathcal{T}'$ annule
ainsi $[\mathcal{T}',[\mathcal{T}',\mathcal{T}']]$, donc
se factorise par une alg\`ebre de Lie nilpotente, donc r\'esoluble. La
restriction de $\rho$ \`a $\mathcal{T}_3$, donc \`a $\mathcal{T}' \simeq
\mathcal{T}_3 / Z(\mathcal{T}_3)$, est semi-simple (cf. \cite{ASSOC} 5.2,
lemme 6). On en d\'eduit que la restriction de $\rho$ \`a $\mathcal{T}'$
est commutative (cf. \cite{BOURB} ch. I \S 5 no. 3 p. 65, cor. 1 du th. 1),
et en particulier $\rho([t_{12},t_{23}]) = 0$. Enfin, si $\k = \Q_l$
on a naturellement $(i) \Rightarrow (iv)$, et la d\'emonstration
de $(iv) \Rightarrow (ii)$ est analogue
\`a celle de $(i) \Rightarrow (ii)$ car
pour tout $\Phi \in \Ass_{\la}(\k)$ il existe $\sigma \in
\Gal(\bar{\Q}|\Q(\mu_{l^{\infty}}))$ tel que $\cc(\sigma. \Phi) \neq
\cc(\Phi)$ d'apr\`es le lemme \ref{ccgal} et la proposition \ref{k3nonnul}.
\end{proof}

Lorsque $\k$ est un corps topologique et la repr\'esentation GT-rigide
consid\'er\'ee est agr\'egeante,
la repr\'esentation projective de $GT_1(\k)$ se d\'ecompose en
caract\`eres continus $\chi : GT_1(\k) \to \k((h))^{\times}$
qui ont les propri\'et\'es suivantes : a) si $L$ est un surcorps de $\k$,
$\chi$ s'\'etend en $\tilde{\chi} : GT_1(L) \to L((h))^{\times}$;
b) pour toute topologie sur $\k$ compatible \`a sa structure de corps,
$\tilde{\chi}$ est continue d\`es que l'inclusion $\k \subset L$
est continue. Nous aurons besoin du lemme suivant.
\begin{lemme} \label{descente}
Soit $\chi : GT_1(\Q) \to \Q((h))^{\times}$ un caract\`ere
v\'erifiant les conditions a) et b), et $g_0$ l'unique \'el\'ement
de $GT_1(\C)$ tel que $g_0 . \Phi_{KZ} = \bar{\Phi}_{KZ}$.
Si $\chi(g_0) \neq 1$, alors $\chi$ est non triviale.
\end{lemme}
\begin{proof}
Pour tout nombre premier $l$, on peut choisir un rationnel $\beta_l
> 0$ tel que $\sqrt{- \beta_l} \in \Q_l$ (par exemple
$\beta_2 = 7$ et $\beta_l = l-1$ si $l \neq 2$). Soit
$L = \Q(\sqrt{-\beta_l}) \subset \C$. Pour la topologie naturelle
de $\C$, les inclusions $\Q \subset L \subset \C$ sont des inclusions
de corps topologiques. Comme $\beta_l >0$, $L$ est dense dans $\C$ donc
$GT_1(L)$ est dense dans $GT_1(\C)$ d'apr\`es le
th\'eor\`eme \ref{theodensite} et il existe donc $g \in GT_1(L)$
tel que $\chi(g') \neq 1$ par continuit\'e. En particulier,
pour tout corps $L'$ isomorphe \`a $L$ il existe $g' \in GT_1(L')$ tel
que $\chi(g') \neq 1$. Prenant pour $L'$ l'extension de $\Q$ dans $\Q_l$
engendr\'ee par $\sqrt{-\beta_l}$ on en d\'eduit qu'existe
$g' \in GT_1(\Q_l)$ tel que $\chi(g') \neq 1$. Comme $\Q$ est dense
dans $\Q_l$ on d\'eduit du th\'eor\`eme \ref{theodensite} qu'existe
$g'' \in GT_1(\Q)$ tel que $\chi(g'') \neq 1$.
\end{proof}

\subsection{Exemples}

Les exemples les plus simples de repr\'esentations GT-rigides et
agr\'egeantes proviennent de l'alg\`ebre d'Iwahori-Hecke de type A,
et sont \'etudi\'ees dans la partie suivante. Dans \cite{DRIN},
th. A (voir aussi \cite{CHARI} prop. 16.1.6), Drinfeld montre
que toutes les composantes (absolument) irr\'eductibles des
repr\'esentations de type Yang-Baxter de $B_n$ sont GT-rigides. Enfin,
il n'est pas difficile de montrer, en reprenant les arguments
de \cite{QUOT} prop. 4, et en utilisant le th. 4 de \cite{QUOT},
que les repr\'esentations de l'alg\`ebre de Birman-Wenzl-Murakami
sont GT-rigides et agr\'egeantes (voir \'egalement \cite{IRRED} th. 2).
Nous d\'eterminons en appendice
la forme g\'en\'erale d'un \'el\'ement de $GT_1(\k)$ modulo des termes
d'ordres au moins 6. Chacun de ces \'el\'ements est congru
\`a l'un des \'el\'ements $g_{a,b}$, d\'efinis dans l'appendice,
pour $a,b \in \k$. Cela nous permet de calculer jusqu'en
degr\'e 5 les valeurs des caract\`eres associ\'es \`a une
telle
repr\'esentation, si l'on sait d\'ecrire une repr\'esentation
infinit\'esimale correspondante. Comme exemple, nous
consid\'erons la repr\'esentation irr\'eductible de
dimension 3 la plus g\'en\'erale de $\mathfrak{B}_3(\k)$, d\'efinie par
$s_1 = \diag(1,1,-1)$,
$t_{12} = \diag(\frac{3+v}{2}, \frac{3-v}{2}, 0)$ et
$$
s_2 = \frac{1}{4v} \left( \begin{array}{ccc}
v-3 & 3(v-1) & 3(v+1) \\
3(v+1) & 3+v & -3(v+1) \\
2v & 2 v \frac{1-v}{1+v} & 2v \\
\end{array} \right)
$$
On a $t_{13} + t_{23} = \diag( \frac{3-v}{2}, \frac{3+v}{2},3)$,
donc cette repr\'esentation irr\'eductible est agr\'egeante si
$v \not\in \{ -3,0,3 \}$. 
Elle est \'egalement GT-rigide, car elle correspond \`a une
repr\'esentation de l'alg\`ebre de Birman-Wenzl-Murakami, ou encore
\`a une repr\'esentation de l'alg\`ebre de Hecke cyclotomique
du groupe de r\'eflexions complexes $G_4$ (cf. \cite{ASSOC}).
Soit $\Phi_0 \in \Ass_1$ et
$R = \widehat{\Phi}(\rho)$. On en d\'eduit deux caract\`eres $Q_{R,2}$
et $Q_{R,3}$ tels que
$$
Q_{R,2}(g_{a,b}) = 1 - \frac{1}{8} av(v^2-9) h^3 - \frac{9}{64} v (v^2-9)
(v^2+7)b  h^5 + \dots
$$
et $
Q_{R,3}(g_{a,b})$ vaut
$$
1 + \frac{1}{16} a(v+9)(v^2-9) h^3 + \frac{9}{128}
(v+5) (v^2-9)(v^2-4v+27)b
(v^2+7) h^5 + \dots
$$

\section{Alg\`ebre de Hecke et caract\`eres de $GT_1(\k)$}

\subsection{Alg\`ebre d'Iwahori-Hecke de type A}
Soit $q \in K^{\times}$ un scalaire non nul. L'alg\`ebre d'Iwahori-Hecke
de type A est not\'ee $H_n(q)$. C'est le quotient de $K B_n$
par la relation $(\sigma_1 -q)(\sigma_1 + q^{-1}) = 0$. On a
$H_n(1) = K \SN$, et $H_n(q)$ est isomorphe \`a $K \SN$ pour $q$ non
racine de l'unit\'e. A toute partition $\alpha$ de $n$ on associe
classiquement, de fa\c con uniforme en $q$, une classe
d'isomorphisme de repr\'esentations de $H_n(q)$, de sorte que la
partition $[n]$ corresponde \`a la repr\'esentation $\sigma_r \mapsto
q \in GL_1(K)$ de $B_n$.

Nos conventions sur les diagrammes de Young sont telles que la partition
$[3,2]$ est repr\'esent\'e par le diagramme \`a deux colonnes
$\Smallruban{ & \ \\ & \\ & \\}$. Dans le tableau de Young standard
$\Smallruban{ 2 & \ \\ 1 & 3 \\}$ la case contenant un 2 est en
colonne 1 et ligne 2. A tout tableau de Young standard $T$ de taille $n$
est associ\'e un $(n-1)$-uplet appel\'e son contenu, d\'efini
comme $(l_i(T)-c_i(T))_{i=2..n}$ o\`u $l_i(T)$ (resp. $c_i(T)$)
d\'esigne la ligne (resp. la colonne) o\`u se trouve
$i$. Cette correspondance est injective, et identifie
les tableaux standard de taille $n$ \`a certains \'el\'ements de $\Z^{n-1}$.
On munit $\Z^{n-1}$ de l'ordre inverse de l'ordre lexicographique,
et l'ensemble des tableaux standard de m\^eme taille de l'ordre
induit.

On d\'ecrit maintenant des expressions matricielles de ces
re\-pr\'e\-sen\-ta\-tions. Soit $\alpha$ une partition de $n$ ; on introduit
le $K$-espace vectoriel de base les tableaux de
Young standard de forme $\alpha$. Pour $r \in [1,n-1]$ on fait agir
$\sigma_r \in B_n$ sur un tableau standard $T$ comme suit. Si $r$ et $r+1$ se
trouvent sur la m\^eme colonne (resp. ligne) de $T$, $\sigma_r.T = qT$
(resp. $\sigma_r.T = -q^{-1} T$). Sinon, soit $T'$ le tableau
(standard) d\'eduit de $T$ par transposition 
de $r$ et $r+1$ dans $T$. Quitte \`a \'echanger les r\^oles de $T$ et $T'$,
on peut supposer que la colonne de $T$ o\`u se trouve $r$ pr\'ec\`ede
celle o\`u se trouve $r+1$, ce qui revient \`a demander $T < T'$. On
note $d$ la \emph{distance axiale}
entre $r$ et $r+1$ dans $T$. Si $r$ (resp. $r+1$) se trouve dans $T$
en colonne $c$ (resp. $c'$) et ligne $l$ (resp. $l'$), elle est d\'efinie
par $d = l-l'+c'-c$, et est donc n\'ecessairement positive. 
On demande que l'action de $\sigma_r$ laisse le plan engendr\'e par
$T$ et $T'$ stable, et que sa restriction soit donn\'ee sur la base
$(T,T')$ par une matrice $M_d^q$ ne d\'ependant que de $d$ et $q$.  On
appelle \emph{mod\`ele matriciel} de $H_n(q)$ toute collection
$(M_d^q)_{d \geq 2}$ telle que la construction pr\'ec\'edente
fournisse une repr\'esentation de $H_n(q)$ de classe correspondant
\`a la partition $\alpha$, pour toute partition $\alpha$ et
tout $n\geq 2$, ceci pour presque tout $q$. Il est imm\'ediat
que chaque $M_d^q$ doit avoir pour trace $q -q^{-1}$ et d\'eterminant
$-1$, et donc \^etre de la forme
$$
\left( \begin{array}{cc} a_d & b_d \\ \frac{1+a_d a'_d}{b_d} & a'_d
\end{array} \right)
$$
avec $a_d \in K$, $b_d \in K^{\times}$, $a'_d = q-q^{-1} - a_d$, et
$1 + a_d a'_d \neq 0$ pour $q$ g\'en\'erique. On v\'erifie
facilement que toute collection $(M_d^q)_{d \geq 2}$ d\'efinie par de
tels
coefficients $a_d$, $b_d$, d\'efinit bien un mod\`ele
matriciel si et seulement si $a_{d+1}(q-a_d) = a_d q^{-1}$ pour tout
$d \geq 2$ avec $a_2 = -1/q(q^2+1)$. On en d\'eduit les
formules g\'en\'erales $a_d = (1-q^2)/q(q^{2d}-1)$,
$a'_d = q^{2d-1}(q^2-1)/(q^{2d}-1)$. Un mod\`ele est donc
arbitrairement d\'etermin\'e par le choix d'une suite
$(b_d)_{d \geq 2}$ d'\'el\'ements de $K^{\times}$. Quel que soit
le choix de ces scalaires, chaque tableau standard est un vecteur
propre commun des $\delta_r, r\in [2,n]$. Si le contenu de $T$
est $(c_2,\dots,c_n)$, $\delta_r.T = q^{2 c_r} T$.

De fa\c con analogue, on appelle \emph{mod\`ele matriciel} du groupe
sym\'etrique une suite $(M_d)_{d \geq 2}$ avec $M_d \in GL_2(\k)$
telle que la construction pr\'ec\'edente donne, pour $q =1$, les
repr\'esentations du groupe sym\'etrique. Deux d'entre eux nous seront
utiles, le mod\`ele
semi-normal
et le mod\`ele orthonormal de Young, respectivement donn\'es par
$$
\frac{1}{d} \left( \begin{array}{cc} -1 & d+1 \\ d-1 & 1 \end{array}
\right) \ \mbox{ et } \ \frac{1}{d} \left( \begin{array}{cc}
-1 & \sqrt{d^2-1} \\ \sqrt{d^2-1} & 1 \end{array} \right)
$$

\subsection{Repr\'esentations infinit\'esimales}

Soit $\rho_0 : \SN \to GL_N(\k)$ une re\-pr\'e\-sen\-ta\-tion du groupe
sy\-m\'e\-tri\-que. Elle s'\'etend en une re\-pr\'e\-sen\-ta\-tion $\rho$ de
$\mathfrak{B}_n(\k)$
par la formule $\rho(t_{ij}) = 2 \rho((i \ j))$. Pour tout $\Phi \in
\Ass_1(\k)$, $\widehat{\Phi}(\rho)(\sigma_1) = \rho(s_1) \exp (h \rho(s_1))$
est semi-simple \`a valeurs propres $q = e^h$ et $-q^{-1}$, donc
$\widehat{\Phi}(\rho)$ se factorise par $H_n(q)$ pour $q = e^h$. Si
$\rho_0$ est absolument irr\'eductible
il en est de m\^eme de $\rho$ et de $\widehat{\Phi}(\rho)$. De plus,
$\rho$ est alors agr\'egeante et $\widehat{\Phi}(\rho)$ permet de
d\'efinir des caract\`eres de $GT_1(\k)$.

Fixons un mod\`ele matriciel $(M_d^1)_{d \geq 2}$ du groupe sym\'etrique,
et \'etendons l'action de $\SN$ sur les tableaux standard en une action
de $\AN$ suivant la formule pr\'ec\'edente. Si $T$ est un tableau standard
de contenu $(c_2,\dots,c_n)$ on a $Y_r.T = 2 c_r T$ pour tout
$r \in [2,n]$. En particulier, si $T$ est un tableau standard contenant
$r$ et $r+1$ dans la m\^eme colonne (resp. la m\^eme ligne), la
droite engendr\'ee par $T$ est stable par $Y_r$ et $t_{r,r+1}$. On fixe
d\'esormais $\Phi \in \Ass_1(\k)$ et on fait agir ($K$-lin\'eairement) $B_n$
sur les tableaux standard selon $\widehat{\Phi}(\rho)$. Comme $\Phi$
est l'exponentielle d'une s\'erie de Lie $\Psi$ sans terme lin\'eaire,
$\Phi(Y_r,t_{r,r+1})$
laisse alors $T$ invariant donc $\sigma_r. T = s_r \exp (h s_r). T$
vaut $q T$ (resp. $-q^{-1} T$). Dans le cas contraire, soit $T'$
le tableau standard d\'eduit de $T$ par transposition de $r$ et $r+1$,
notons $(c'_2,\dots,c'_n)$ son contenu et supposons $T < T'$. Le plan
engendr\'e par $T$ et $T'$ est stable par $s_r$, $Y_r$ et $Y_{r+1}$,
donc par $\sigma_r$. L'expression de $s_r \exp (hs_r)$ dans
la base $(T,T')$ vaut $M_d^1 \exp ( h M_d^1)$, avec $d$ la distance
axiale associ\'ee au couple $(T,T')$. D'autre part, $Y_r + Y_{r+1}$
commute \`a $Y_r$ et $t_{r,r+1}$, et $Y_r = (Y_r + Y_{r+1})/2 +
(Y_r - Y_{r+1})/2$, donc $\Phi(Y_r,t_{r,r+1}) = \Phi((Y_r - Y_{r+1})/2,t_{r,
r+1})$ car $\Phi = \exp \, \Psi$ avec $\Psi \in \mathcal{L}(\k)$ sans
terme lin\'eaire. Comme
$t_{r,r+1}$ agit par $2 M_d^1$ et $(Y_r-Y_{r+1})/2$ par la matrice
diagonale $\eta_d$ de coefficients $(c_r-c'_r,c_{r+1}-c'_{r+1}) = (d,-d)$,
l'expression de $\sigma_r$ dans la base $(T,T')$ ne d\'epend que de
$d$ et $q$, et fournit ainsi un mod\`ele matriciel $(M_d^q)_{d \geq 2}$
de $H_n(q)$, tel que $M_d^q \equiv M_d^1$ modulo $h$. En particulier,
\`a tout $\Phi \in \Ass_1(\k)$ est ainsi attach\'e un mod\`ele
$(M_d^q(\Phi))_{d \geq 2}$ et une suite de scalaires $b_d(\Phi) \in
K^{\times}$, $d \geq 2$. On remarque que, comme $(M_d^1)^2=1$, pour tout $\Phi \in \Ass_1(\k)$, on
a
$
2M_d^q(\Phi) = (q-q^{-1}) + (q+q^{-1}) Q M_d^1 Q^{-1}
$
avec $Q = \Phi( 2 h M_d^1,h \eta_d)$. Comme nous n'utiliserons que deux
mod\`eles du groupe sym\'etrique, nous noterons $b_d^s(\Phi)$
et $b_d^o(\Phi)$ les scalaires $b_d(\Phi)$ correspondant respectivement
aux mod\`eles semi-normal et orthogonal. Comme ces deux mod\`eles
sont conjugu\'es par la matrice diagonale de coefficients $(\sqrt{d+1},
\sqrt{d-1})$, on en d\'eduit que ces deux nombres sont li\'es par
la relation $b_d^s(\Phi)\sqrt{d-1} = b_d^o(\Phi) \sqrt{d+1}$. Un calcul 
explicite montre d'autre part ais\'ement que
$$
\frac{d}{d+1} b_d^s(\Phi) = 1 + \frac{h^2}{6} + 16 \cc(\Phi) d h^3 + o(h^3).
$$ 

\subsection{L'associateur KZ}

Dans sa th\`ese \cite{JORGE}, Jorge Gonz\'alez-Lorca a calcul\'e
$M_d^q(\Phi_{KZ})$, donc $b_d(\Phi_{KZ})$. Pour la commodit\'e du
lecteur mais aussi parce que nous aurons \'egalement besoin de
conna\^\i tre $b_d(\bar{\Phi}_{KZ})$, nous reprenons ici les derni\`eres
\'etapes de ce calcul. On a ici $\k = \C$, $A = \C[[h]]$, $K= \C((h))$.
On note $\hh = h/(2\ii \pi)$ et l'on suppose $M_d^1$ donn\'e par le
mod\`ele semi-normal. On introduit la fonction de trois variables
complexes
$$
\Gamma_3(a,b,c) = \frac{\Gamma(1-2a) \Gamma(1+2b)}{\Gamma(1+b-a+\delta)
\Gamma(1+b-a-\delta)}$$
o\`u $\delta$ d\'esigne une racine carr\'ee de $a^2+b^2+2c$ et
$\Gamma$ est la fonction Gamma d'Euler. Comme $\Gamma(1+z)$ est analytique
au voisinage de 0, on peut d\'efinir $\Gamma_3(\alpha,\beta,\gamma) \in K$
pour $\alpha,\beta,\gamma \in h A$ et $F(a,b) = \Gamma_3(\hh a, \hh b, -2
\hh^2) \in K$ pour $a,b \in \k$. Pour tout $X \in M_N(\k)$ semi-simple et $x \in Sp(X)$
on note $P_{x,X} \in M_N(\k) \subset M_N(K)$ le projecteur sur l'espace propre
$\Ker(X-x)$ naturellement associ\'e et on d\'efinit $F(U,V)$,
pour $U,V \in M_N(\k)$ semi-simples, comme la somme des termes $F(u,v) P_{u,U}
P_{v,V}$ sur tous les couples $(u,v) \in Sp(U) \times Sp(V)$. Pour traiter
simultan\'ement les cas de $\Phi_{KZ}$ et $\bar{\Phi}_{KZ}$, on introduit
un param\`etre $\eps = \pm 1$ et l'on pose $F_{\eps}(a,b) =
F(\eps a, \eps b)$.

Pour all\'eger les notations, notons
$U = \eta_d$, $V = 2 M_d^1$ et $S = M_d^1$. En utilisant le
fait que $UV+VU = -4$ et les m\'ethodes de Riemann et Kummer sur
l'\'equation hyperg\'eom\'etrique, J. Gonz\'alez-Lorca montre
dans \cite{JORGE} que $\Phi_{KZ}(\tilde{U},\tilde{V}) = F_{\eps}(V,U)$,
$\Phi_{KZ}(\tilde{V},\tilde{U}) = F_{\eps}(U,V)$,
o\`u l'on note $\tilde{U} = \eps h U$, $\tilde{V} = \eps h V$. On en
d\'eduit $M_d^q(\Phi)$ pour $\Phi = \Phi_{KZ}$ (si $\eps = 1$)
et $\Phi = \bar{\Phi}_{KZ}$ (si $\eps = -1$) comme suit. Soit $P_{\pm}
= P_{\pm 1,S} = P_{\pm 2,V}$. On a
$$
\begin{array}{lcl}
F_{\eps}(U,V) & = & F_{\eps}(U,2)P_+ + F_{\eps}(U,-2)P_- \\
F_{\eps}(V,U) & = & P_+ F_{\eps}(2,U) + P_- F_{\eps}(-2,U) \\
\end{array}
$$
donc $\Phi_{KZ}(\tilde{V},\tilde{U})S \Phi_{KZ}(\tilde{U},
\tilde{V})$ est \'egal \`a la diff\'erence
$$
\left(F_{\eps}(U,2)P_+ F_{\eps}(2,U) \right) - \left( 
 F_{\eps}(U,-2) P_- F_{\eps}(-2,U) \right)
$$
et ainsi son coefficient
en premi\`ere ligne et deuxi\`eme colonne vaut
$$
\frac{d+1}{2d} \left[ F_{\eps}(d,2) F_{\eps}(2,-d) + F_{\eps}(d,-2)
F_{\eps}(-2,-d) \right].
$$
On en d\'eduit
$$
\begin{array}{lcl}
b_d^s(\Phi_{KZ}) & = & \frac{d+1}{2d} \frac{q+q^{-1}}{2} ( Z(d,2) + Z(d,-2) ) \\
b_d^s(\bar{\Phi}_{KZ}) & = & \frac{d+1}{2d} \frac{q+q^{-1}}{2} ( Z(-d,-2) + Z(-d,2) ) 
\end{array}
$$
avec $Z(a,b) = F(a,b) F(b,-a)$, donc $b_d^s(\Phi_{KZ})/
b_d^s(\bar{\Phi}_{KZ})$ vaut $\tilde{J}(d+1)\tilde{J}(d-1)/\tilde{J}(d)^2$
avec $\tilde{J}(x) = J(2 \hh x)$ et
$$
J(x) = \frac{\Gamma(1+x)}{\Gamma(1-x)} = \exp \left(
-2 \gamma x - 2 \sum_{n=1}^{\infty} \frac{\zeta(2n+1)}{2n+1} x^{2n+1} \right)
$$
o\`u $\gamma$ est la constante d'Euler. Ainsi,
$$
(*) \ \ \ \frac{b_d^s(\Phi_{KZ})}{b_d^s(\bar{\Phi}_{KZ})} =
\exp \left( -2 \sum_{n=1}^{\infty} \frac{\zeta(2n+1)}{2n+1} 2^{2n+1}
\hh^{2n+1} Q_n(d) \right)
$$
avec $Q_n(d) = (d+1)^{2n+1} + (d-1)^{2n+1} - 2d^{2n+1} $.

\subsection{Associateurs pairs}

Soit $\Phi \in \Ass_1^0(\k)$, et choisissons pour $M_d^1$ le mod\`ele
orthonormal de Young, de telle sorte que $M_d^1$ est simultan\'ement
orthogonale et sym\'etrique pour tout $d \geq 2$. On en d\'eduit
que $M_d^1 \exp(h M_d^1)$ est \'egalement sym\'etrique. Soit $\Psi$
la s\'erie de Lie sans terme lin\'eaire telle que $\Phi = \exp \, \Psi$.
Comme $\eta_d$ est sym\'etrique et $x \mapsto - ^t x$ est un automorphisme
de l'alg\`ebre de Lie $\mathfrak{gl}_N(K)$, la transpos\'ee de $\Psi(h \eta_d,
h M_d^1)$ vaut $-\Psi(-h \eta_d,-h M_d^1)$. On en d\'eduit
$$
^t \Phi(h \eta_d, h M_d^1) = \Phi(-h \eta_d,-h M_d^1)^{-1} =
\Phi(h \eta_d, h M_d^1)^{-1}
$$
et qu'ainsi $M_d^q(\Phi)$ est sym\'etrique. On a donc $b_d^o(\Phi)^2 =
1 + a_d a'_d$. Comme $1+a_d a'_d$ est congru \`a $(d^2-1)/d \neq 0$
modulo $h$, il existe un unique $b_d^0(\Phi) \in A$ v\'erifiant
cette \'equation qui soit de plus congru \`a $\sqrt{d^2-1}/d$ modulo $h$.
Avec les notations des $q$-analogues $[n]_q = (q^n - q^{-n})/(q-q^{-1})$,
on trouve
$
b_d^o(\Phi) = \sqrt{ [d+1]_q [d-1]_q} /[d]_q
$
d'o\`u
$$
b_d^s(\Phi) = \frac{\sqrt{d+1} \sqrt{ [d+1]_q [d-1]_q}}{\sqrt{d-1} [d]_q}
= \frac{d+1}{d} \left( 1 + \frac{h^2}{3!} - \frac{4d^2-1}{5!} h^4 + \dots \right)
$$
Remarquons qu'avec ces notations, $a_d = -q^{-d}/[d]_q$, $a'_d = q^d/[d]_q$.
Comme $\k \subset \R$ on note comme cons\'equence imm\'ediate de ces calculs
\begin{prop} \label{modunitaire} Un mod\`ele matriciel de $H_n(q)$ unitaire au sens de
\cite{ASSOC} est donn\'e par
$$
M_d^q(\Phi) = \frac{1}{[d]_q} \left( \begin{array}{cc}
-q^{-d} & \sqrt{ [d+1]_q [d-1]_q } \\
\sqrt{ [d+1]_q [d-1]_q } & q^d \end{array} \right)
$$
\end{prop}
Afin de comparer $b_d^s(\Phi)$ et $b_d^s(\Phi_{KZ})$, on
introduit les fonctions
$$
I(z) = \Gamma(1+z) \Gamma(1-z) = \frac{\pi z}{\sin(\pi z)} \ \ \ \tilde{I}
(z) = I(2z \hh).
$$
Comme $(q^n - q^{-n})/2 = \ii \sin(2n\pi \hh) = nh/\tilde{I}(n)$, on a
$$
b_d^s(\Phi) = \frac{d+1}{d} \frac{\tilde{I}(d)}{\sqrt{\tilde{I}(d+1)
\tilde{I}(d-1)}}
$$
d'o\`u, en utilisant $(q+q^{-1}) \tilde{I}(2)/2 \tilde{I}(1) = 1$
et l'expression de $\Gamma(1 \pm x)$ en fonction de $I$ et $J$,
on obtient
$$
\frac{b_d^s(\Phi_{KZ})}{b_d^s(\Phi)}  = \tilde{H}(d) \mbox{ et }
\frac{b_d^s(\overline{\Phi}_{KZ})}{b_d^s(\Phi)}  = \frac{1}{\tilde{H}(d)}
$$
avec 
$$
\tilde{H}(d) = \frac{\sqrt{\tilde{J}(d-1) \tilde{J}(d+1)}}{\tilde{J}(d)}
= \exp \left( -2 \sum_{n=1}^{\infty} \frac{\zeta(2n+1)}{2n+1} 2^{2n+1}
\hh^{2n+1} Q_n(d) \right)
$$
et $Q_n(d) = (d+1)^{2n+1} + (d-1)^{2n+1} - 2d^{2n+1} $.

Les formules obtenues nous disent en particulier que le logarithme
de $b_d^s(\Phi_{KZ}) / b_d^s(\Phi)$ admet pour d\'eveloppent
limit\'e
$$  \frac{-2 \ii \zeta(3) d}{\pi^3}
h^3 + \frac{2 \ii d (2d^2+1) \zeta(5)}{\pi^5} h^5 - \frac{2\ii d (1+3d^4+5d^2)
\zeta(7)}{\pi^7} h^7 + o(h^8)
$$

\subsection{Caract\`eres}

Choisissons un mod\`ele matriciel de $H_n(q)$ associ\'e \`a une suite
$(b_d)_{d \geq 2} = (b_d(\Phi))_{d \geq 2}$ pour un certain $\Phi \in
\Ass_{\la}(\k)$. Pour toute partition $\alpha$ de $n$, la
repr\'esentation irr\'eductible $R$ associ\'ee sur les tableaux
standard de forme $\alpha$ est GT-rigide et $R(\delta_2),\dots,R(\delta_n)$
engendre les matrices diagonales, comme on le d\'eduit imm\'ediatement
de l'expression de $\delta_r$ pour $2 \leq r \leq n$ et de l'invariance de
$\sigma_1$ sous l'action de $GT_1(\k)$. En particulier, prenons pour $R$
la repr\'esentation correspondant \`a la partition $[2,1^{n-2}]$ -- c'est
la repr\'esentation de Burau (r\'eduite). On note $e_1,\dots,e_{n-1}$
sa base suivant l'ordre d\'efini sur les tableaux standard,
\begin{figure}
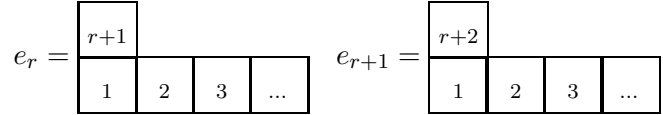

$$
e_r = \ruban{ r+1 \\ 1 & 2 & 3 & \dots \\} \ \ \ 
e_{r+1} = \ruban{ r+2 \\ 1 & 2 & 3 & \dots \\}
$$ 
\caption{Base de la repr\'esentation de Burau}
\end{figure}
et $\chi_d = Q_{R,d}$ pour $2 \leq d \leq n-1$. Pour tout $g \in GT_1(\k)$,
on prend pour repr\'esentant de $Q_R(g)$ dans $GL_N(K)$
la matrice diagonale de coefficients $(1,\chi_2(g),\chi_2\chi_3(g),\dots)$.
Notant par abus $\sigma'_r$ (resp. $\sigma_r$) l'action de $\sigma_r \in
B_n$ suivant $g.R$ (resp. $R$) on a, sur le plan $(e_r,e_{r+1})$,
$$
\left( \begin{array}{cc} * & b'_{r+1} \\
* & *  \end{array} \right) = \sigma'_{r+1} =
\lmm \chi(g) & 0 \\ 0 & \chi'(g)
\rmm ^{-1} 
\lmm * & b_{r+1} \\ * & * \rmm 
\lmm \chi(g) & 0 \\ 0 & \chi'(g)
\rmm
$$
avec $\chi = \chi_2\dots\chi_r$ et $\chi' = \chi_2 \dots \chi_{r+1}$,
d'o\`u $b'_{r+1} = b_{r+1} \chi_{r+1}(g)$. Comme d\'emontr\'e en 3.2.3, les
$\chi_r$ sont \`a valeurs dans $A^{\times} = \k[[h]]^{\times}$.
Supposons maintenant
$\Phi \in \Ass_1(\k)$ et $b_d = b_d^s(\Phi)$. Si $\Phi' = g.\Phi$
on a alors $b_d^s(\Phi') = b_d^s(\Phi) \chi_d(g)$.
En particulier, si $g \in GT_1(\C)$ est tel que $\bar{\Phi}_{KZ} =
g.\Phi_{KZ}$, $\chi_d(g)^{-1}$ est donn\'e par la formule $(*)$. Comme,
\`a $d \geq 2$ fix\'e, $Q_n(d) \sim 2n(2n+1)d^{2n-1}$ pour $n$ grand,
les fonctions de $n$ d\'efinies sont lin\'eairement
ind\'ependantes et on d\'eduit imm\'ediatement du lemme \ref{descente} :
\begin{prop} Les caract\`eres $\chi_d : GT_1(\Q) \to (\Q[[h]])^{\times}$
pour $d \geq 2$ sont alg\'ebriquement ind\'ependants
sur $\Q$ en tant que fonctions $GT_1(\Q) \to \Q[[h]]$.
\end{prop}
Pour conclure la d\'emonstration du th\'eor\`eme C,
soit $\alpha$ une
partition quelconque de $n$, $R$ la repr\'esentation de $H_n(q)$ 
correspondante suivant le mod\`ele associ\'e \`a $b_d^s(\Phi)$
et $N$ sa dimension. Si $T$ est un tableau standard de taille $n$,
et $i<j$ avec $c_i(T) > c_j(T)$, on note $d_T(i,j) = l_j(T)-l_i(T)+c_i(T)
-c_j(T) > 0$. Un repr\'esentant dans $GL_N(K)$ de $Q_R(g)$ pour
$g \in GT_1(\k)$ est alors donn\'e par la matrice diagonale
qui \`a $T$ de forme $\alpha$ associe $\chi_T(g) T$ avec
$$
\chi_T(g) = \prod_{\stackrel{i<j}{c_i(T) > c_j(T)}} \chi_{d_T(i,j)}(g).
$$
En effet, il suffit de montrer que cette matrice conjugue les
re\-pr\'e\-sen\-ta\-tions associ\'ees \`a $b_d^s(\Phi)$ et $b_d^s(\Phi')$,
c'est-\`a-dire les expressions de $\sigma_r$ pour $1 \leq r \leq n-1$
suivant ces deux repr\'esentations. Fixant un tel $r$, il suffit
de le montrer sur les plans $(T_1,T_2)$ o\`u $T_2$ est d\'eduit de
$T_1$ par l'\'echange de $r$ et $r+1$. Or, si $T_1 <T_2$ et en notant
$d$ la distance axiale associ\'ee, on a $d = d_{T_2}(r,r+1)$
et
$\chi_{T_2}(g) = \chi_{T_1}(g) \chi_d(g)$, donc $M_d^q(\Phi)$
est conjugu\'ee par la matrice diagonale de coefficients $(1,\chi_d(g))$,
ce qui donne bien $M_d^q(\Phi')$.

La relation entre les caract\`eres issus de deux diagrammes sym\'etriques
l'un de l'autre est la suivante. Soit $T$ un tableau de forme
$\alpha$, $T'$ son sym\'etrique de forme $\alpha'$. Pour $d <0$,
notons par convention $\chi_d = 1$. Alors
$$
\chi_T \chi_{T'} = \prod_{\stackrel{i<j}{c_i(T) > c_j(T)}}
\chi_{d_T(i,j)} \prod_{\stackrel{i<j}{l_i(T) > l_j(T)}}
\chi_{d_T(j,i)}
$$ est encore \'egal au produit des $\chi_{d_T(i,j)}
\chi_{d_T(j,i)}$ sur tous les couples $i<j$ tels que $c_i(T) > c_j(T)$
ou (exclusif) $l_i(T) > l_j(T)$. De tels couples sont en bijection
avec les crochets du diagramme $\alpha$, c'est-\`a-dire les couples
$(x,y)$ de cases de $\alpha$ o\`u $x$ (resp. $y$) est en colonne $c$
(resp. $c'$) et ligne $l$ (resp. $l'$) avec $c<c'$ et $l>l'$ :
\`a $(x,y)$ on associe le couple $(i,j)$ o\`u $i$ (resp. $j$) est
le minimum (resp. le maximum) de leurs contenus dans $T$ ; inversement,
la condition sur $i$ et $j$ signifie que les deux cases qui les
contiennent forment un crochet. Pour un tel crochet $(x,y)$ on
d\'efinit sa longueur $\delta = \delta(x,y) = c'-c+l-l'$. Alors
$$
\chi_{d_T(i,j)} \chi_{d_T(j,i)} = \chi_{\delta} \chi_{-\delta} = \chi_{
\delta}
$$
et
$$
\chi_T \chi_{T'} = \prod_{\mbox{crochets }(x,y)} \chi_{\delta(x,y)}
$$
d\'epend seulement de $\alpha$, et plus de $T$.

Soit $\varepsilon$ l'automorphisme involutif de $K$ d\'efini par
$f(h) \mapsto f(-h)$ et $U_N^{\varepsilon}(K) = \{ x \in GL_N(K) \ \mid \ 
x^{-1} = ^t \varepsilon(x) \}$ le groupe unitaire associ\'e. Si l'on
part d'une repr\'esentation orthogonale
$\rho_0 : \mathfrak{S}_n \to O_N(\k)$, les repr\'esentations
$\widehat{\Phi}(\rho)$ de l'alg\`ebre de Hecke obtenues se
factorisent par $U_N^{\varepsilon}(K)$ (cf. \cite{ASSOC} section 3.2.2).
On en d\'eduit que les caract\`eres $\chi_d$ se factorisent
par $U_1^{\varepsilon}(A) = \{ x \in A \ \mid \ \varepsilon(x) = x^{-1} \}$,
et que les $\chi_d(g)$ sont des exponentielles de s\'eries
\emph{impaires} en $h$ (sans terme lin\'eaire d'apr\`es le calcul de $b_d^s(\Phi)$
\`a l'ordre 3 en 4.2).
 
\subsection{R\'esonances}

Remarquons que si une partition de $n$ est donn\'ee, les caract\`eres
$\chi_T$ pour $R$ parcourant les tableaux standard de forme $\alpha$
ne sont pas n\'ecessairement distincts, comme le montre
l'exemple de $\alpha = [3,2]$ : un repr\'esentant de $Q_R(g)$ dans
$GL_5(K)$ est alors donn\'e par la matrice diagonale de
coefficients $(1,\chi_3(g),\chi_2\chi_3(g),\chi_2\chi_3(g),\chi_2^2 \chi_3(g))$.

\begin{prop} \label{equerres} La repr\'esentation de $H_n(q)$ associ\'ee \`a un
diagramme $\alpha$ est sans r\'esonances si et seulement si $\alpha$
est une \'equerre ou \'egal \`a $[2,2]$.
\end{prop}
\begin{proof}
Soit $R$ la repr\'esentation de Burau (r\'eduite), et $e_1,\dots,
e_{n-1}$ la base choisie pr\'ec\'edemment. Comme les diagrammes en \'equerres
correspondent aux puissances ext\'erieures de $R$, il faut montrer que
$\Lambda^r R$ est sans r\'esonances pour $r \in [0,n-1]$. Fixons $r$ et
associons \`a toute
partie $I = \{i_1,\dots,i_r \} \subset [1,n-1]$ de cardinal $r$ avec
$i_1<\dots<i_r$
le vecteur $e_I = e_{i_1} \wedge \dots \wedge e_{i_r}$. L'ensemble des
vecteurs de ce type forme une base de $\Lambda^r R$. Un
rel\`evement lin\'eaire de $Q_R$ est donn\'e par l'action
lin\'eaire $g.e_i = 
\psi_i(g) e_i$, avec $\psi_1 = 1$, $\psi_2 = \chi_2$, \dots,
$\psi_{n-1} = \chi_2 \dots \chi_{n-1}$. Un rel\`evement dans
$GL(\Lambda^r R)$ de l'action projective de $GT_1(\k)$ sur $\Lambda^r R$
est donc donn\'e par $g \mapsto Q(g)$ avec
$$
Q(g) e_I = \left( \prod_{i \in I} \psi_i \right) (g) \, e_I
= \left( \prod_{d=2}^{n-1} \chi_d^{f_d(I)} \right) (g) \, e_I
$$ 
o\`u $f_d(I) = \# \{ i \in I \ \mid \ i \geq d \}$. Comme les
$\chi_d$ sont alg\'ebriquement ind\'ependants, il s'agit donc de
montrer que, si $f_d(I) = f_d(J)$ pour $d \in [2,n-1]$, alors
$I= J$, au moins si $\# I = \#J = r$. Mais dans ce cas, on a \'egalement
$\# I = f_1(I) = f_1(J) = \# J$, d'o\`u l'on d\'eduit imm\'ediatement
$I = J$, puisque $i \in I \Leftrightarrow f_i(d) > f_{i-1}(d)$.

Le cas du diagramme $[2,2]$ se ram\`ene \`a celui de l'\'equerre $[2,1]$,
puisque ce deuxi\`eme est la restriction \`a $B_3 \subset B_4$ du
premier.

D'apr\`es la relation entre $\chi_T$ et $\chi_{T'}$ pour deux tableaux
$T$ et $T'$ sym\'etriques l'un de l'autre, la repr\'esentation associ\'ee \`a $\alpha$ sera
sans r\'esonances si et seulement si celle associ\'ee \`a $\alpha'$
l'est. Il reste donc \`a montrer que, si $\alpha \geq [3,2]$,
il existe deux tableaux distincts $T_1$ et $T_2$ de forme $\alpha$ tel que
$\chi_{T_1} = \chi_{T_2}$. Consid\'erons les deux tableaux de
la figure \ref{figtabexc} : on remplit d'abord les cases du
sous-diagramme $[3,2]$ comme indiqu\'e avec les nombres de 1 \`a 5, puis
l'on remplit les autres cases dans l'ordre standard avec les nombres
restant. On a alors
$$
\begin{array}{lcccl}
\chi_{T_1} & = & (\chi_3 \dots \chi_{n_1 -1})(\chi_2 \dots \chi_{n_1-3})
& = & \chi_2 \chi_3^2 \dots \chi_{n_1 -3}^2 \chi_{n_1-2} \\
& = &  (\chi_2 \dots \chi_{n_1-2})(\chi_3 \dots \chi_{n_1-3}) & = & 
\chi_{T_2} \\
\end{array}
$$

\begin{figure}
$$
T_1 = \ruban{ n_1 & \ & \ & \  \\ \vdots & n_2 & \ & \ \\
6 & \vdots & \ & \ \\ 5 & 1+n_1 & n_3 & \ \\2 & 4 & \vdots & \ \\
1 & 3 & 1+n_2 & \dots \\} \ \ 
T_2 = \ruban{ n_1 & \ & \ & \  \\ \vdots & n_2 & \ & \ \\
6 & \vdots & \ & \ \\ 4 & 1+n_1 & n_3 & \ \\3 & 5 & \vdots & \ \\
1 & 2 & 1+n_2 & \dots \\} \ \ 
$$
\caption{\label{figtabexc}
Tableaux tels que $\chi_{T_1} = \chi_{T_2}$}
\end{figure}

\end{proof}

\section{Appendice : Cocycles de $GT_1(\k)$}

On fait agir $\alpha \in \GM(\k) = \k^{\times}$ sur
$f(h) \in K = \k((h))$ par $(f.\alpha)(h) = f(\alpha h)$.
Cette action laisse stable l'anneau de valuation discr\`ete
$A = \k[[h]]$, son id\'eal maximal $A_0 = hA$ et
$A_1 = 1 + h A = \exp A_0$. Le morphisme naturel $\pi : GT_1(\k)
\to \kt$ de noyau $GT_1(\k)$ induit ainsi une action de $GT(\k)$
sur $K$, et donc sur $PGL_N(K)$. De mani\`ere \`a
obtenir, non plus seulement des repr\'esentations projectives de
$GT_1(\k)$ mais des 1-cocycles (a priori non ab\'eliens) de $GT(\k)$ \`a
valeurs dans $PGL_N(K)$, on peut renforcer la d\'efinition 3.1 comme suit :
\begin{defi}
Une repr\'esentation $R \in V_N(\k)$ est dite \emph{fortement}
GT-rigide si $R$ est absolument irr\'eductible et, pour tout
$g \in GT(\k)$, $g.R$ est isomorphe \`a $\sigma \mapsto R(\sigma).\pi(g)$.
\end{defi}

Si $R \in V_N(\k)$ est fortement GT-rigide, $R$ est en particulier
GT-rigide. De plus il existe alors,
pour tout $g \in GT(\k)$, un \'el\'ement $Q_R(g) \in PGL_N(K)$
tel que
$$
\forall \sigma \in B_n \ \ \ (g.R)(\sigma) = Q_R(g)^{-1}
\left( R(\sigma).\pi(g) \right) Q_R(g).
$$
On v\'erifie facilement que $Q_R(g_1g_2) = (Q_R(g_1).\pi(g_2)) Q_R(g_2)$,
c'est-\`a-dire que $Q_R \in Z^1(GT(\k),PGL_N(K))$. Si de plus
$R$ est agr\'egeante, on en d\'eduit des cocycles $\chi \in Z^1(GT(\k),
A_1)$ ainsi que des cocycles
$$
\psi = \log \chi \in Z^1(GT(\k),A_0).
$$
Ces cocycles permettent enfin de d\'efinir des repr\'esentations
$GT(\k) \to GL_2(A)$
par
$$
g \mapsto \left( \begin{array}{cc} \pi(g) & \psi(g) \\ 0 & 1 \end{array}
\right)
$$
Parmi les repr\'esentations $R \in V_N(\k)$ qui sont fortement
GT-rigides, on trouve notamment les repr\'esentations irr\'eductibles
de l'alg\`ebre d'Iwahori-Hecke de type A et celles de l'alg\`ebre de
Birman-Wenzl-Murakami.

Si la restriction de $\chi$ \`a $GT_1(\k)$ est non triviale (cf. prop.
\ref{trivia}), $\chi$ n'est pas un cobord et la repr\'esentation
de $GT(\k)$ dans $GL_2(A)$ qui s'en d\'eduit n'est pas scind\'ee.
De fa\c con g\'en\'erale, on montre inversement que les \'el\'ements non
nuls de $H^1(GT(\k), A_1)$ proviennent n\'ecessairement
de caract\`eres non triviaux de $GT_1(\k)$.

\begin{prop} Le morphisme de restriction $H^1(GT(\k),A_1) \to
\Hom(GT_1(\k),A_1)$ est injectif.
\end{prop}
En effet, la suite exacte d'inflation-restriction en cohomologie
des groupes associ\'ee \`a la suite exacte $1 \to GT_1(\k) \to GT(\k)
\to \kt \to 1$ est
$$
0 \to H^1(\kt,A_1) \to H^1(GT(\k),A_1) \to H^1(GT_1(\k),A_1)
$$
avec $H^1(GT_1(\k),A_1) = \Hom(GT_1(\k),A_1)$. Il suffit
donc de montrer que $H^1(\kt,A_1) = 0$. Pour ce faire notons, pour
$m \geq 0$, $\k_m$ le $\kt$-module $\k$ pour l'action $(\alpha,
x) \mapsto \alpha^m x$.
\begin{lemme}
Si $m \geq 1$, $H^1(\kt,\k_m) = 0$.
\end{lemme}
\begin{proof}
Si $f \in Z^1(\kt,\k_m)$ et $\alpha_1,\alpha_2 \in \k$ sont
tels que $\alpha_1^m \neq 1$ et $\alpha_2^m \neq 1$,
$$
f(\alpha_1 \alpha_2) = \alpha_1^m f(\alpha_2) + f(\alpha_1)
= f(\alpha_2 \alpha_1) = \alpha_2^m f(\alpha_1) + f(\alpha_2)
$$
donc $f(\alpha)/(\alpha^m -1)$ ne d\'epend pas du choix de
$\alpha$, pourvu que $\alpha^m \neq 1$. Comme $\mathrm{car.} \k = 0$,
$\kt \neq \mu_m(\k)$ et il existe $x \in \k$ bien d\'efini
tel que
$f(\alpha) = (\alpha^m -1)x = \alpha.x -x$
pour tout $\alpha \in \kt \setminus \mu_m(\k)$. Il reste \`a montrer
que $f(\alpha) = 0$ pour tout $\alpha \in \mu_m(\k)$. Or la
restriction de $f$ \`a $\mu_m(\k)$ est dans $\Hom(\mu_m(\k),
\k) = 0$ puisque $\k$ est sans torsion, donc $f$ est un
cobord et $H^1(\kt,\k_m) = 0$.
\end{proof}

On en d\'eduit
$$
H^1(\kt,A_1) \simeq H^1(\kt,A_0) \simeq \prod_{m\geq 1} H^1(\kt,\k_m) = 0,
$$
le premier isomorphisme d\'ecoulant de l'isomorphisme de $\kt$-modules
entre $(A_0,+)$ et $(A_1,\times)$ donn\'e par l'exponentielle,
et le second du fait que $A_0 \simeq \prod_{m \geq 1} \k_m$
en tant que $\kt$-module.

\section{Appendice : Calculs en degr\'e 5}

Pour faire les calculs explicitement, on introduit
des bases form\'ees de mon\^omes de Lie des composantes homog\`enes
de $\mathcal{A}(\k)$ jusqu'en degr\'e 5. D'autre part, on
note $a \equiv b$ si $\omega(a-b) \geq 6$.
\begin{table}
$$
\begin{array}{|lclrclrcl|}
\hline
w_3 & = & [x,y] & w_4 & = & [x,[x,y]] & w_5 & = & [y,[y,x]] \\
w_6 & = & [x,[x,[x,y]]] & w_7 & = & [y,[x,[x,y]]] & w_8 & = &
[y,[y,[x,y]]] \\
w_9 & = & [x,[x,[x,[x,y]]]] & w_{10} & = & [y,[x,[x,[x,y]]]] &
w_{11} & = & [y,[y,[x,[x,y]]]] \\
w_{12} & = & [y,[y,[y,[x,y]]]] & w_{13} & = & [[x,y],[x,[x,y]]] &
w_{14} & = & [[x,y],[y,[x,y]]] \\
\hline
\end{array}
$$
\caption{Bases de $\mathcal{L}(\k)$ en degr\'e au plus 5}
\end{table}
Notons $\psi_1$ et $\psi_2$ les \'el\'ements de $\mathfrak{grt}_1(\k)$ correspondant
aux \'el\'ements $f_1$ et $f_2$ de l'interpr\'etation \og hamiltonienne
\fg\ de Drinfeld dans \cite{DRIN}. On exprime facilement $\psi_1 = w_5 - w_4$, et
$$
\frac{\partial f_2}{\partial x} = -2 w_{10} -2w_{11} - 2 w_{12} + 2 w_{13}
+ w_{14}.
$$
L'\'el\'ement $\psi_2$ se d\'eduit de $\partial f_2 / \partial x$
par la transformation involutive $F(x,y) \mapsto F(x,-x-y)$,
dont la matrice dans la base $w_9,\dots,w_{14}$ des \'el\'ements
de degr\'e 5 s'\'ecrit
$$
\left( \begin{array}{cccccc}
-1 & 1 & -1 & 1 & 0 & 0 \\
0 & 1 & -2 & 3 & 0 & 0 \\
0 & 0 & -1 &3 & 0 & 0 \\
0 & 0 & 0 & 1 & 0 & 0 \\
0 & 0 & -1 & 2 & 1 & -1 \\
0 & 0 & 0 & 1 & 0 & -1 \\
\end{array}
\right)
$$
d'o\`u
$$
\psi_2 = -2 w_9 -4w_{10} - 4 w_{11} - 2w_{12} -w_{13} -3w_{14}.
$$
On pose $\psi_{a,b} = a\psi_1 +b \psi_2 \in \mathfrak{grt}_1(\k)$ pour $a,b \in\k$,
et on note $\Psi_{a,b} \in GRT_1(\k)$ son exponentielle au sens du groupe.
Si $\Phi \in \Ass_1(\k)$, on sait que l'on peut \'ecrire
$\Phi \equiv 1 + \frac{1}{24} w_3 - \cc \psi_1 + \phi_4 + \phi_5$
avec $\phi_4,\phi_5$ homog\`enes de degr\'es respectifs 4 et 5. Alors
$$
\Phi . \Psi_{a,b} = \exp(s_{\psi_{a,b}})(\Phi)
\equiv \Phi + a \psi_1 + b \psi_2 + \frac{a}{24}(w_3 \psi_1 + [[\psi_1,x],
y]),
$$
et de plus $[[\psi_1,x],y]] = -w_{11} - w_{10}$. D'autre part
la formule de Le et Murakami \cite{LEMURAK} dit
$$
\Phi_{KZ} \equiv 1 + \frac{w_3}{24} + \tilde{\zeta}(3) \psi_1
+ \phi_4^{KZ} + \frac{1}{2} \tilde{\zeta}(5) \psi_2 -
\frac{1}{24} \tilde{\zeta}(3) (w_{10}+w_{11} + w_3\psi_1)
$$
avec $\tilde{\zeta}(n) = \zeta(n)/(2\ii \pi)^n$,
et $\phi_4^{KZ} = (-w_7 - 4w_6 - 4w_8 + 5 w_3^2)/5760$.
On en d\'eduit d'une part que 
$$
\Phi_0 = \Phi_{KZ}.\Psi_{-\tilde{\zeta}
(3),-\tilde{\zeta}(5)/2} \equiv 1 + \frac{1}{24} w_3 + \phi_4^{KZ}
$$
est pair jusqu'en degr\'e 6, et d'autre part que tout associateur
a les m\^emes termes de degr\'e 2 et 4. 

Pour obtenir l'expression g\'en\'erale d'un \'el\'ement de $GT_1(\k)$
jusqu'en degr\'e 5, il suffit de calculer $g_{a,b} = \iota_{\Phi_0}
(\Psi_{a,b})$, uniquement d\'etermin\'e par $g_{a,b}. \Phi_0 = \Phi_0
. \Psi_{a,b}$. Posons $g_{a,b} = \exp F$. A des termes d'ordre
au moins 4 pr\`es, on a $\Phi_0 = 1 + w_3/24$, $\Phi_0 e^x \Phi_0^{-1}
= e^x + [w_3,x]/24$, et $\log \Phi_0 e^x \Phi_0^{-1} = x - w_4/24$.
Notons $G_4$, $G_5$ les composantes homog\`enes de degr\'e 4 et 5
de $g_{a,b}$. On a vu que $g_{a,b} \equiv 1 + a \psi_1 + G_4 + G_5$,
d'o\`u $F \equiv a \psi_1 + G_4 + G_5$. Or $\psi_1 = w_5-w_4$,
et
$$ \left\lbrace \begin{array}{lcl}
w_4(x-w_4/24,y)& \equiv & w_4 + (w_{10} +2w_{13})/24 \\
w_5(x-w_4/24,y)& \equiv & w_5 - w_{11}/24 \\
\end{array} \right.
$$
d'o\`u
$$
F(\log \Phi_0 e^x \Phi_0^{-1}),y) \equiv a \psi_1 + G_4 + G_5 - \frac{a}{24}
(w_{10}+ w_{11} +2 w_{13}).
$$
D'autre part $\Phi_0 . \Psi_{a,b} \equiv \Phi_0 + a \psi_1 + b \psi_2 +
a (w_3 \psi_1 - w_{11} - w_{10})$ d'o\`u
$$
\begin{array}{lcl}
(\Phi_0. \Psi_{a,b}). \Phi_0^{-1} 
& = & 1 + a \psi_1 + b \psi_2 + a([w_3,\psi_1] - w_{11} - w_{10})/24 \\
& = & 1 + a \psi_1 + b \psi_2 - a(w_{10} + w_{11} + w_{13} + w_{14})/24 \\
\end{array}
$$
et son logarithme vaut
$$
a \psi_1 + b \psi_2 - a(w_{10} + w_{11} +
w_{13} + w_{14})/24 = F(\log \Phi_0 e^x \Phi_0^{-1}),y)
$$
d'o\`u $G_4 = 0$, et $G_5 = b \psi_2 + a (w_{13}- w_{14})/24$ soit
$$
G_5 = - \left( 2b w_9 + 4b w_{10}
+ 4b w_{11} + 2b w_{12} + (b - \frac{a}{24}) w_{13} + (b + \frac{a}{24}) w_{14} \right)
$$
On peut d'autre part calculer
$$
\chi_d(g_{a,b}) = \frac{b^s_d(\Phi_0.\Psi_{a,b})}{b^s_d(\Phi_0)} = 1 - 16ad h^3 -128db(1+2d^2)h^5 + \dots
$$

Pour faire la comparaison avec les termes $\kappa_m^*(\sigma)$
d'Ihara, il faut calculer, pour chaque mon\^ome de Lie $w_r$, les
\'el\'ements $w'_r = \pi_{ab} \circ p_x ( w_r(\log(1+x),\log(1+y)))$.
On trouve
$$
\left\lbrace
\begin{array}{lcl}
w'_4 & \equiv & \frac{1}{2} y^2x^2-\frac{1}{3}y^3x^2-x^2y+x^3y-\frac{1}{2}
y^2x^3-\frac{11}{12}x^4y \\
w'_5 & \equiv & \frac{1}{3}y^2x^3+y^2x-\frac{1}{2}y^2x^2+\frac{1}{2}y^3x^2
-y^3x+\frac{11}{12}y^4x \\
\end{array}
\right.
$$
et les images de $w_9,\dots,w_{14}$ sont respectivement $-x^4 y, -y^2 x^3,
-y^3 x^2$, $- y^4 x,0,0$.
On d\'eduit alors de la formule g\'en\'erale de $g_{a,b}$
et de
$$
\begin{array}{rl}
\psi^{\sigma}_{ab} \equiv &
1 + \frac{\kappa_3^*(\sigma)}{2} (yx^2 + y^2x - yx^3 - y^2x^2 - y^3 x) \\
& + ( \frac{1}{24} \kappa_5^*(\sigma) + \frac{11}{24} \kappa_3^*(\sigma))(x^4y + y^4 x) \\
& + (\frac{5}{12} \kappa_3^*(\sigma) + \frac{1}{12} \kappa_5^*(\sigma))(y^3x^2 + y^2 x^3) \\
\end{array}
$$
que $a = \kappa_3^*(\sigma)/2$ et $b = \kappa_5^*(\sigma)/48$.

\bigskip

Si $l=3$, $\kappa_3^*(\sigma) = \kappa_3(\sigma)/8$,
$\kappa_5^*(\sigma) = \kappa_5(\sigma)/80$. On peut v\'erifier
$$
\frac{1}{24} (\kappa^*_5 + 11 \kappa^*_3) = 
\frac{1}{8.80} \frac{1}{3}(\kappa_5 + 110 \kappa_3) \in \Z_3
$$
parce que modulo 3 $\kappa_5 + 110 \kappa_3 \equiv 111\kappa_3 \equiv 0$,
$$
\frac{1}{12} (5\kappa^*_3 +  \kappa^*_5) = 
\frac{1}{4.80} \frac{1}{3}(50\kappa_3 + \kappa_5) \in \Z_3
$$
car $50\kappa_3 + \kappa_5 \equiv 51 \kappa_3 \equiv 0$ modulo 3.

Si $l = 2$, on a $(\kappa_5^* + 11 \kappa_3^*)/24 = (3
\kappa_5 + 341 \kappa_3)/(2^3.279) \in \Z_2$
car $3
\kappa_5 + 341 \kappa_3 \equiv 8 \kappa_3 \equiv 0$ modulo 8.
De m\^eme $(5\kappa_3^* + \kappa_5^*)/12 = (155 \kappa_3 + 3 \kappa_5)/(2^2.
279) \in \Z_2$ parce que $155 \kappa_3 + 3 \kappa_5 \equiv 3 (\kappa_3
+ \kappa_5 \equiv 3.2 \kappa_3 \equiv 0$ modulo 4. On a donc v\'erifi\'e
que les coefficients de $\psi^{ab}(\sigma)$ appartiennent \`a $\Z_l$.
On d\'eduit \'egalement de ces formules qu'au moins les premiers
coefficients de $\chi_d(\sigma)$, $-16ad = -8d\kappa_3^*(\sigma)$
et $-128bd(1+2d^2) = -8 \kappa_5^*(\sigma)d(1+2d^2)/3$ appartiennent
\`a $\Z_l$, puisque $d(1+2d^2)$ est divisible par 3 pour tout $d \in \N$.

\end{document}